# Induced Saturation of Graphs

Maria Axenovich[a] and Mónika Csikós [*a]

[a]Institute of Algebra and Geometry, Karlsruhe Institute of Technology

**Abstract**

A graph $G$ is $H$-saturated for a graph $H$, if $G$ does not contain a copy of $H$ but adding any new edge to $G$ results in such a copy. An $H$-saturated graph on a given number of vertices always exists and the properties of such graphs, for example their highest density, have been studied intensively.

A graph $G$ is $H$-induced-saturated if $G$ does not have an induced subgraph isomorphic to $H$, but adding an edge to $G$ from its complement or deleting an edge from $G$ results in an induced copy of $H$. It is not immediate anymore that $H$-induced-saturated graphs exist. In fact, Martin and Smith (2012) showed that there is no $P_4$-induced-saturated graph. Behrens *et al.* (2016) proved that if $H$ belongs to a few simple classes of graphs such as a class of odd cycles of length at least 5, stars of size at least 2, or matchings of size at least 2, then there is an $H$-induced-saturated graph.

This paper addresses the existence question for $H$-induced-saturated graphs. It is shown that Cartesian products of cliques are $H$-induced-saturated graphs for $H$ in several infinite families, including large families of trees. A complete characterization of all connected graphs $H$ for which a Cartesian product of two cliques is an $H$-induced-saturated graph is given. Finally, several results on induced saturation for prime graphs and families of graphs are provided.

# 1 Introduction

The notion of induced saturation is relatively new. A more studied subgraph saturation problem does not impose the "induced" condition. A graph $G$ is called $H$-*saturated* if $G$ has no subgraph isomorphic to $H$ but adding any edge to $G$ from its complement creates a subgraph isomorphic to $H$. The problem of finding dense $H$-saturated graphs has a long history dating back to 1907, when Mantel [8] determined the maximum number of edges in a triangle-free graph of a fixed order. In 1941, Turán [11] determined maximum number of edges in a $K_k$-saturated graph on $n$ vertices, which gave rise to the classical extremal function $\text{ex}(n, H)$, that is the largest number of edges in an $H$-saturated graph on $n$ vertices. The study of $H$-saturated graphs with smallest possible number of edges was carried out later by Erdős, Hajnal and Moon [4], who determined the minimum number of edges in a $K_k$-saturated graph on $n$ vertices. Since then, graph and hypergraph saturation became an attractive area of research, see the extensive dynamic survey of Faudree, Faudree and Schmitt [5]. Note that $H$-saturated graphs always exist and can be constructed greedily from an empty graph by adding edges as long as no copy of $H$ appears.

---

[*]Corresponding author
*E-mail addresses:* maria.aksenovich@kit.edu and monika.csikos@kit.edu.

A graph $H$ is an *induced subgraph* of a graph $G$ if $V(H) \subseteq V(G)$ and for any $x, y \in V(H)$, $xy \in E(H)$ if and only if $xy \in E(G)$. We say that $G$ contains an *induced copy* of $H$ and write $H \subseteq_I G$, if $G$ contains an induced subgraph isomorphic to $H$.

The notion of induced extremal [10] and induced saturation [9] numbers were first introduced for generalized graphs, called trigraphs, that allow for "flexible" edges which can be added or omitted when one looks for an induced copy of the desired graph. In particular, the induced saturation and induced extremal numbers were defined as the minimum and maximum number of "flexible" edges in an $H$-induced-saturated trigraph.

In this paper, we work with the following definition of induced saturation. We say that $G$ is *$H$-induced-saturated*, and write $G \hookrightarrow\!\!\!\!\rightarrow H$, if $G$ does not have an induced subgraph isomorphic to $H$, but deleting any edge of $G$ or adding any edge to $G$ from its complement creates an induced copy of $H$. Note that if $G \hookrightarrow\!\!\!\!\rightarrow H$, then $\overline{G} \hookrightarrow\!\!\!\!\rightarrow \overline{H}$.

We use the following notation: $K_n$, $C_n$, and $P_n$ for a complete graph, a cycle, and a path on $n$ vertices, respectively; $K_{n,m}$ for a complete bipartite graph with parts of sizes $n$ and $m$; $C'_n$ for a union of $C_n$ and $K_2$, that share exactly one vertex; $\widetilde{C}_n$ for a union of $C_n$ and an edge $e$ where the endpoints of $e$ are distance 2 apart in $C_n$.

Martin and Smith [9] introduced the notion of induced saturation number and as a by-product they proved that there is no graph $G$ such that $G \hookrightarrow\!\!\!\!\rightarrow P_4$. We provide a short proof of this result in Section 6.

All previously known examples of graphs $H$ for which there exists an $H$-induced-saturated graph were given by Behrens, Erbes, Santana, Yager, and Yaeger [1].

**Theorem 1** (Behrens et al.[1]). *Let $n \geq 2$, $k \geq 3$. If $H$ is one of the following graphs: $C'_3$, $K_{1,n}$, $nK_2$, $C_4$, $C_{2k-1}$, $C'_{2k}$, $\widetilde{C}_{2k}$, then there is a graph $G$ such that $G \hookrightarrow\!\!\!\!\rightarrow H$.*

In this paper, we introduce two large classes of graphs, $\mathcal{H}(n)$ and $\mathcal{F}(n)$, including a family of trees, such that for any graph $H$ in these classes, there is an $H$-induced-saturated graph $G$, which is a Cartesian product of cliques. We need a few definitions to state our results.

The *Cartesian product* of $n$ cliques of order $k$ is denoted by $\overset{n}{\square} K_k$ and is formally defined as a graph whose vertex set is a set of all $n$-tuples with entries from a $k$-element set such that two tuples are adjacent if and only if they differ in exactly one coordinate. We say that a graph is *Hamming* if it is isomorphic to $\square^n K_k$ for some $n, k \in \mathbb{N}$. We write $K_k \square K_k$ for $\square^2 K_k$. We define the *dimension*, $\dim(G)$ of a graph $G$ as the smallest integer $n$ such that $G \subseteq_I \square^n K_k$ holds for some $k \in \mathbb{N}$. If no such $n$ exists, the dimension is defined to be infinity.

We say that a graph is a balanced *p-legged spider* of height $q$ if it is a union of $p$ paths, $L_1, \ldots, L_p$ each of length $q$ such that they share exactly one vertex, $r$, called the *head* of the spider, i.e., $V(L_i) \cap V(L_j) = \{r\}$ for $1 \leq i < j \leq p$. The vertices of degree 1 are called the *leaves* and the paths $L_1, \ldots, L_p$ are called the legs of the spider. A connected subgraph of a balanced spider is a spider. A spider is a $p$-legged spider of height at least $q$ if it contains a balanced $p$-legged spider of height $q$ as a subgraph and has maximum degree $p$.

Now we define the classes $\mathcal{H}(n)$ and $\mathcal{F}(n)$.

- For any $n \geq 3$, let $\mathcal{H}(n)$ be a class of connected graphs that have a cut vertex $r$ of degree $n+1$ so that $H - r$ has $n+1$ connected components, each of dimension at most $n-1$. See Figure 1a.

- For any $n \geq 3$, let $\mathcal{F}(n)$ be the class of graphs $H$ that can be obtained as the union of graphs $F$ and $T$ such that $\dim(F) \leq n-1$, $T$ is a balanced $(n+1)$-legged spider of height $n+2$, and $|V(F) \cap V(T)| = V'$, where $V'$ is the set of leaves of $T$. See Figure 1b.

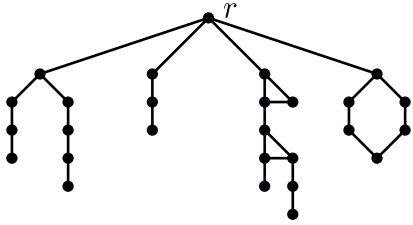
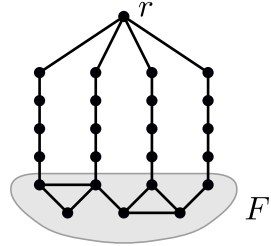

(a) A graph from $\mathcal{H}(3)$.  (b) A graph from $\mathcal{F}(3)$.

Figure 1: Examples of graphs from classes $\mathcal{H}(n)$ and $\mathcal{F}(n)$.

**Theorem A.** *For any $n \geq 3$ and $H \in \mathcal{H}(n)$, there is $k \in \mathbb{N}$ such that $\square^n K_k \hookrightarrow\!\!\!\!\!\rightarrow H$.*

**Corollary 2.** *For any $n \geq 3$, let $\mathcal{T}(n)$ be the class of trees with a unique vertex of maximum degree $n+1$ such that its removal results in components that have maximum degree at most $n-1$. Then for any $n \geq 3$ and $T \in \mathcal{T}(n)$, there is $k \in \mathbb{N}$ so that $\square^n K_k \hookrightarrow\!\!\!\!\!\rightarrow T$.*

**Theorem B.** *For any $n \geq 3$ and $H \in \mathcal{F}(n)$, there is $k \in \mathbb{N}$ such that $\square^n K_k \hookrightarrow\!\!\!\!\!\rightarrow H$.*

We say that $H'$ is obtained from $H$ by a *k-subdivision* if $V(H') = V(H) \dot\cup \left(\dot\bigcup_{e \in E(H)} V_e\right)$, for any $xy \in E(H)$, $|V_{xy}| = k$, $V_{xy} \cup \{x, y\}$ induces a path with endpoints $x$ and $y$, and there are no other edges in $E(H')$. Informally, $H'$ is the graph obtained from $H$ by subdividing each edge with $k$ new vertices.

**Corollary 3.** *Let $n \geq 4$ and let $H$ be a graph with a unique vertex $v$ of maximum degree $n + 1$ such that $H - v$ has maximum degree at most $n - 1$. Let $H'$ be the graph obtained from $H$ by an $(n+1)$-subdivision. Then there is $k \in \mathbb{N}$ such that $\overset{n}{\square} K_k \hookrightarrow\!\!\!\!\twoheadrightarrow H'$.*

Note that the same result can be easily proved when $H'$ is obtained from $H$ by subdividing each edge with at least $n + 1$ vertices.

In Section 5, we give a characterization of connected graphs $H$ for which there exists an integer $k$ so that $K_k \square K_k \hookrightarrow\!\!\!\!\twoheadrightarrow H$.

Our results give infinite classes of graphs $H$ for which $H$-induced-saturated graphs exist. These classes include some easily describable families of graphs that complement the list of Theorem 1:

**Proposition 4.** *Let $k \geq 3$ and let $H$ be one of the following graphs: a spider with $k$ legs, $\widetilde{C}_{2k+1}$, or the union of a path $P$ and a cycle $C$ such that $P$ has length at least 1 and shares its endpoint and no other vertex with $C$. Then there is a graph $G$ such that $G \hookrightarrow\!\!\!\!\twoheadrightarrow H$.*

The paper is organized as follows. In Section 2, we give some basic properties of Hamming graphs. In Sections 3 and 4, we prove Theorems A and B, respectively. In Section 5, we describe those connected graphs $H$ for which there is an $H$-induced-saturated graph that is a Cartesian product of two cliques of equal size. We also prove Proposition 4 in Section 5. In the last two sections, we give some general observations about induced saturation of prime graphs and state some concluding remarks.

## 2 Preliminaries and properties of Hamming graphs

Let $G$ and $H$ be graphs. An injective mapping $\sigma : V(H) \to V(G)$ is an *induced embedding* of $H$ into $G$ if $uv \in E(H)$ holds if and only if $\sigma(u)\sigma(v) \in E(G)$. Note that $G$ contains an induced copy of $H$ if and only if $H$ has an induced embedding into $G$.

Let $G$ be a connected graph and $u, v \in V(G)$. The *distance* between $u$ and $v$ in $G$, denoted by $\operatorname{dist}_G(u, v)$ is the length of a shortest $u, v$-path in $G$. We shall denote the maximum degree of a graph $G$ by $\Delta(G)$ and its chromatic index by $\chi'(G)$.

The *Hamming distance* between two vectors $\mathbf{u}, \mathbf{v} \in \mathbb{R}^n$, denoted by $\operatorname{dist}_{\text{Ham}}(\mathbf{u}, \mathbf{v})$, is the number of coordinates where $\mathbf{u}$ and $\mathbf{v}$ differ.

Let **0** denote a vector with all zero entries and $\mathbf{e}_i$ denote a vector with all entries zero, except for the $i^{\text{th}}$ which is 1. We shall denote the vertices of $\square^n K_k$, that are elements of a set $X^n$, $|X| = k$ by bold letters and the $i^{th}$ coordinate of a vector $\mathbf{v}$ by $[\mathbf{v}]_i$. Let $[n] = \{1, 2, \ldots, n\}$ and $[n, m] = \{k \in \mathbb{N} : n \leq k \leq m\}$.

**Definition 5.** *Let $G$ be a graph, $n \in \mathbb{N}$. We call a function $c : E(G) \to X$ a **nice coloring** with $n$ colors if $|X| = n$,*

  i) *all triangles are monochromatic, and*

  ii) *for any two distinct, nonadjacent vertices $u, v$ of $G$ there exist two different colors $i$ and $j$ such that $i$ and $j$ appear on each induced $u, v$-path.*

The next observation follows from the second property of a nice coloring and the fact that for adjacent vertices $u, v$, the edge $uv$ is the only induced $u, v$-path.

**Observation 6.** *Let $G$ be a graph and $c$ be a nice coloring of $E(G)$. Then for any two distinct vertices $u, v$ of $G$, there exists a color $i$ such that $i$ appears on each induced $u, v$-path.*

The class of graphs with finite dimension was characterized by Klavžar and Peterin [7]:

**Theorem 7** (Coloring criterium [7]). *Let $G$ be a graph and $n \in \mathbb{N}$. Then $G$ has dimension at most $n$ if and only if $G$ has a nice coloring with $n$ colors.*

Note that the smallest graph with infinite dimension is $K_4 - e$ and that if $H \subseteq_I G$, then $\dim(H) \leq \dim(G)$.

**Lemma 8.** *Let $n \geq 2$ and $F$ be a forest with $\Delta(F) = n$. Then $\dim(F) = \chi'(F) = n$, in particular, $F \not\subseteq_I \square^{n-1} K_k$ for any $k \in \mathbb{N}$.*

*Proof.* Let $F$ be a forest with $\Delta(F) = n$. By the coloring criterium (Theorem 7), it is enough to show that $F$ has a nice coloring with $n$ colors and cannot be nicely colored with a smaller number of colors. Observe that an edge-coloring of a forest is nice if and only if it is proper. Hence $\dim(F) = \chi'(F)$. It is easy to see that $F$ does not have a proper edge-coloring with less than $n$ colors since the $n$ edges incident to a vertex of maximum degree need to receive different colors. On the other hand, $E(F)$ can be greedily properly colored with $n$ colors. Thus $\dim(F) = \chi'(F) = n$. □

The following lemma is an easy consequence of the fact that Hamming graphs are edge-transitive and distance-transitive.

**Lemma 9.** *Let $k, n \geq 2$ and $G = \square^n K_k$. Then*

- *for any $e, f \in E(G)$, the graphs $G - e$ and $G - f$ are isomorphic;*

- for any $e, f \in E(\overline{G})$ with $e = \mathbf{uv}$, $f = \mathbf{xy}$ such that $\text{dist}(\mathbf{u}, \mathbf{v}) = \text{dist}(\mathbf{x}, \mathbf{y})$, the graphs $G + e$ and $G + f$ are isomorphic.

**Lemma 10.** *Let $n, k \in \mathbb{N}$, $G = \overset{n}{\square} K_k$, and $H$ be a graph. Then $G \hookrightarrow\!\!\!\!\rightarrow H$ if and only if*

*S1)* $H \not\subseteq_I G$,

*S2) for any edge $e \in E(\overline{G})$, $H \subseteq_I G + e$, and*

*S3) there is and edge $f \in E(\overline{H})$ such that $H + f \subseteq_I G$.*

*Proof.* By definition, $G \hookrightarrow\!\!\!\!\rightarrow H$ if and only if
1) $H \not\subseteq_I G$,
2) for any $e \in E(\overline{G})$, $H \subseteq_I G + e$, and
3) for any $e \in E(G)$, $H \subseteq_I G - e$.
Thus it is enough to prove that 3) is equivalent to S3).

Assume that $H \not\subseteq_I G$ and 3) holds. Let $e = uv \in E(G)$ and let $\sigma : V(H) \to V(G)$ be an induced embedding of $H$ into $G - e$. Since $H \not\subseteq_I G$, the vertex set $\sigma(V(H))$ has to contain both $u$ and $v$. Since $\sigma$ is an induced embedding, $\sigma^{-1}(u)\sigma^{-1}(v) \notin E(H)$. Let $f = \sigma^{-1}(u)\sigma^{-1}(v)$, then $\sigma$ is an induced embedding of $H + f$ into $G$. Thus S3) holds.

Let S3) hold. Let $f = u'v' \notin E(H)$ so that $H + f \subseteq_I G$ and let $\sigma : V(H) \to V(G)$ be an induced embedding of $H + f$ into $G$. Since $\sigma$ is an induced embedding, $\sigma(u')\sigma(v') \in E(G)$. Let $e = \sigma(u')\sigma(v')$, then $\sigma$ is an induced embedding of $H$ into $G - e$. Thus $H \subseteq_I G - e$. By Lemma 9, $H \subseteq_I G - e'$ for any $e' \in E(G)$ and so 3) holds. □

**Lemma 11.** *Let $k, n \in \mathbb{N}$, $G = \overset{n}{\square} K_k$, and $H$ be a graph with $H \subseteq_I G$. Then for any $u \in V(H)$ and for any $\mathbf{v} \in V(G)$, there exists an induced embedding $\sigma : V(H) \to V(G)$ so that $\sigma(u) = \mathbf{v}$.*

*Proof.* We assume that $V(G) = \{0, 1, \ldots, k-1\}^n$. Let $\tau : V(H) \to V(G)$ be an arbitrary induced embedding of $H$ into $G$ and let $\mathbf{w} = \tau(u)$. If $\mathbf{w} = \mathbf{v}$, then $\sigma = \tau$ is a desired embedding. If $\mathbf{w} \neq \mathbf{v}$, define $\mu : V(G) \to V(G)$ as $\mu(\mathbf{x}) = \mathbf{x} - \mathbf{w} + \mathbf{v} \pmod{k}$. Then $\mu$ is a bijection that preserves Hamming distance, and thus $\mu$ is an automorphism of $G$. Therefore $\sigma = \mu \circ \tau$ is an induced embedding of $H$ into $G$ and $\sigma(u) = \mathbf{w} - \mathbf{w} + \mathbf{v} = \mathbf{v}$. □

## 3 Proof of Theorem A

Let $H \in \mathcal{H}(n)$, i.e., $H$ has a cut vertex $r$ of degree $n + 1$ such that $H - r$ has $n + 1$ connected components, each of dimension at most $n - 1$. We denote the components of $H - r$ by $H_1, H_2, \ldots, H_{n+1}$ and define $\{v_i\} = V(H_i) \cap N(r)$ for all $i \in [n+1]$, see Figure 2. Let $\ell \in \mathbb{N}$ be the smallest integer so that $H_i \subseteq_I \overset{n-1}{\square} K_\ell$ for all $i \in [n+1]$.

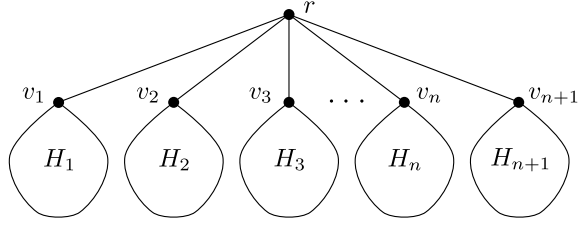

Figure 2: The structure of graphs in $\mathcal{H}(n)$.

Let $G = \square^n K_k$ be the Hamming graph on a vertex set $\{0, 1, \ldots, k-1\}^n$, where $k \in \mathbb{N}$ is to be specified later. We shall show that $G \hookrightarrow\!\!\!\twoheadrightarrow H$. By Lemma 10, $G \hookrightarrow\!\!\!\twoheadrightarrow H$ if and only if

S1) $H \not\subseteq_I G$,

S2) for any edge $e \in E(\overline{G})$, $H \subseteq_I G + e$, and

S3) there is an edge $f \in E(\overline{H})$ such that $H + f \subseteq_I G$.

The statement S1) holds since the vertex set $\{r\} \cup N(r)$ induces a star $K_{1,n+1}$ in $H$, hence by Lemma 8, $H \not\subseteq_I G$.

Now we prove the statement S2). Let $e \in E(\overline{G})$ and $G^+ = G + e$. We need to show that $H \subseteq_I G^+$. Denote the endpoints of $e$ by $\mathbf{u}$ and $\mathbf{v}$, and let $q = \mathrm{dist}_G(\mathbf{u}, \mathbf{v})$, $2 \leq q \leq n$. By Lemma 9, we can suppose that $\mathbf{u} = \mathbf{0}$ and $\mathbf{v}$ has the form $\mathbf{v} = (\underbrace{2, \ldots, 2}_{q}, 0, \ldots, 0)$.

We shall provide an explicit embedding of $H$ into $G^+$ using some auxiliary embeddings $\sigma_i$ of $H_i$ into suitable Hamming graphs $F_i$. We start by defining the graphs $F_1, \ldots, F_{n+1}$. Consider the sets $W_1, \ldots, W_{n+1}, X$ with $W_i \subseteq \{0, 2, 3, \ldots, k-1\}$, $X \subseteq \{2, 3, \ldots, k-1\}$ such that $|W_i| = |X| = \ell$, $W_i \cap W_j = \{0\}$ for all $j \neq i$, $2 \in X \cap W_{n+1}$ and for any $i \in [n]$, $X \cap W_i = \emptyset$. Clearly, such sets $W_i$ and $X$ exist if $k \geq (\ell-1)(n+1) + 3$. Let $F_i$ be the Hamming graph whose vertex set is the set of vectors with $n-1$ components such that if $i \in [n]$, then each coordinate has value from $W_i$ and if $i = n+1$, then the first coordinate is from $X$ and all other coordinates are from $W_{n+1}$. Since $F_i$ is isomorphic to $\square^{n-1} K_\ell$ for any $i \in [n+1]$, there is an induced embedding $\sigma_i$ of $H_i$ into $F_i$. By Lemma 11, we can suppose that $\sigma_i(v_i) = \mathbf{0}$ for all $i \in [n]$ and that $\sigma_{n+1}(v_{n+1}) = (\underbrace{2, \ldots, 2}_{q-1}, 0, \ldots, 0)$. Note that $q - 1 \geq 1$, which implies that the first coordinate of $\sigma_{n+1}(v_{n+1})$ is equal to 2.

Define an embedding $\sigma : V(H) \to V(G^+)$ such that $\sigma(r) = \mathbf{0}$ and for $u \in V(H_i)$, $\sigma(u)$ is obtained by fixing its value in one coordinate and using the values of $\sigma_i(u)$ in the other coordinates. Formally,

(a) if $u \in V(H_1)$, then let $[\sigma(u)]_1 = 1$ and $[\sigma(u)]_j = [\sigma_1(u)]_{j-1}$ for all $2 \leq j \leq n$;

(b) if $u \in V(H_i)$ with $2 \leq i \leq n$, then let $[\sigma(u)]_i = 1$, $[\sigma(u)]_j = [\sigma_i(u)]_j$ for all $1 \leq j < i$, and $[\sigma(u)]_j = [\sigma_i(u)]_{j-1}$ for all $i < j \leq n$;

(c) if $u \in V(H_{n+1})$, then let $[\sigma(u)]_1 = 2$ and $[\sigma(u)]_j = [\sigma_{n+1}(u)]_{j-1}$ for all $2 \leq j \leq n$.

See Figure 3 for an illustration. Note that $\sigma(v_i) = \mathbf{e}_i$ for $i \in [n]$ and $\sigma(v_{n+1}) = \mathbf{v}$.

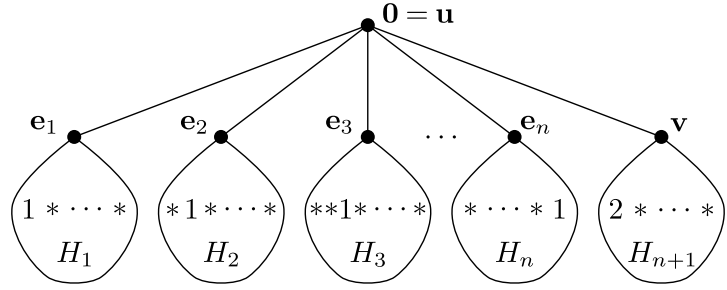

Figure 3: The scheme of the embedding $\sigma$ of $H$ in $G^+$, where $* \in \{0, 2, 3, \ldots, k-1\}$ correspond to the values of $\sigma_i$'s.

**Claim.** The mapping $\sigma$ is an induced embedding of $H$ in $G^+$.

It follows from (a), (b), (c), and the injectivity of the mappings $\sigma_i$ that $\sigma$ is injective. We shall show that for any pair $x, y$ of vertices of $V(H)$, $xy \in E(H)$ if and only if $\sigma(x)\sigma(y) \in E(G^+)$.

By (a), (b), and (c), for any $i \in [n+1]$ and $x, y \in V(H_i)$, $\sigma_i(x)$ and $\sigma_i(y)$ have Hamming distance 1 if and only if $\sigma(x)$ and $\sigma(y)$ have Hamming distance 1. Since the $\sigma_i$'s are induced embeddings, $G^+[\sigma(V(H_i))]$ is isomorphic to $H_i$ for all $i \in [n+1]$.

Also, for any $i \in [n+1]$, $\sigma(r)$ is adjacent to $\sigma(v_i)$ and for any vertex $x$ of $H_i - v_i$, $\sigma(x)$ has at least two nonzero coordinates, hence it is not adjacent to $\sigma(r)$ in $G^+$.

It remains to check that two vertices from distinct $H_i$'s are mapped to two non-adjacent vertices of $G^+$. Consider $x \in V(H_i)$ and $y \in V(H_j)$ with $1 \leq i < j \leq n+1$. Since $\sigma(x)\sigma(y) \neq \mathbf{uv}$, it suffices to show that $\sigma(x)$ and $\sigma(y)$ differ in at least 2 coordinates. If $j \neq n+1$, then by (a) and (b), $[\sigma(x)]_i = 1 \neq [\sigma(y)]_i$ and $[\sigma(y)]_j = 1 \neq [\sigma(x)]_j$ with $i \neq j$. If $j = n+1$ and $2 \leq i \leq n$, then by (b) and (c), $[\sigma(x)]_i = 1 \neq [\sigma(y)]_i$, and $[\sigma(y)]_1 = 2 \neq [\sigma(x)]_1$, where $i \neq 1$. If $j = n+1$ and $i = 1$, then by (a) and (c), $[\sigma(x)]_1 = 1 \neq 2 = [\sigma(y)]_1$ and $[\sigma(x)]_2 \neq [\sigma(y)]_2$ since $[\sigma(x)]_2 \in W_1$, $[\sigma(y)]_2 \in X$, and $W_1 \cap X = \emptyset$.

This proves the Claim and the statement S2).

Finally, we prove the statement S3). Let $f = v_n v_{n+1}$ and $H^+ = H + f$. We shall show that $H^+ \subseteq_I G$, using the coloring criterium (Theorem 7).

Since $H_i \subseteq_I \square^{n-1} K_\ell$, by Theorem 7, there is a nice coloring of $H_i$ with $n-1$ colors. For $i \in [n]$, let $c_i$ be a nice coloring of $H_i$ with colors from $[n] \setminus \{i\}$ and let $c_{n+1}$ be a nice coloring of $H_{n+1}$ with colors from $[n-1]$. Define a coloring $c : E(H^+) \to [n]$ as follows:

$$c(e) = \begin{cases} c_i(e) : & e \in E(H_i), \\ i : & e = rv_i, \ 1 \leq i \leq n, \\ n : & e = rv_{n+1}, \\ n : & e = v_n v_{n+1}. \end{cases}$$

Now we verify that $c$ is a nice coloring. Let $T$ be a triangle in $H^+$. Then either $V(T) = \{r, v_n, v_{n+1}\}$ and so each edge of $T$ has color $n$ or $T$ is contained in $H_i$ for some $i \in [n+1]$ and so $T$ is monochromatic by the niceness of $c_i$.

Let $x, y$ be two distinct, nonadjacent vertices of $H^+$. We need to show that there are two distinct colors that appear on each induced $x, y$-path.

If $x, y \in V(H_i)$ for some $i \in [n+1]$, then all $x, y$-paths are contained in $H_i$ and thus niceness of $c_i$ implies that there is a pair of distinct colors present on all these paths.

If $x \in V(H_i)$, $y \in V(H_j)$, $1 \le i < j \le n+1$, $\{i, j\} \ne \{n, n+1\}$, then any $x, y$-path contains the edges $rv_i$ and $rv_j$, thus contains the distinct colors of these edges.

If $x \in V(H_n)$ and $y \in V(H_{n+1})$, then all induced $x, y$-paths contain $v_n v_{n+1}$, so they contain color $n$. Since $x$ and $y$ are not adjacent in $H^+$ and $v_n v_{n+1} \in E(H^+)$, $x \ne v_n$ or $y \ne v_{n+1}$. Assume, without loss of generality, that $x \ne v_n$. By the niceness of $c_n$ an Observation 6, there is a color $i \in [n-1]$ that appears on every induced $x, v_n$-path, thus all induced $x, y$-paths contain the colors $i$ and $n$.

Finally, if $x = r$ and $y \in V(H_i) \setminus \{v_i\}$ with $i \in [n+1]$, then any induced $x, y$-path contains the edge $rv_i$ and an induced $v_i, y$-path in $H_i$. Hence any induced $x, y$-path contains the color $c(rv_i)$ and since $c_i$ is a nice coloring with colors from $[n] \setminus \{c(rv_i)\}$, by Observation 6, all induced $v_i, y$-paths contain a color $j \in [n] \setminus \{c(rv_i)\}$. Thus all induced $x, y$-paths contain the colors $j$ and $c(rv_i)$.

Therefore $c$ is a nice coloring of $H^+$ with $n$ colors and so by Theorem 7, $\dim(H^+) \le n$, i.e., there exists $m \in \mathbb{N}$ so that $H^+ \subseteq_I \overset{n}{\square} K_m$.

Letting $k = \max\{(\ell - 1)(n+1) + 3, \ m\}$, we see that both statements S2) and S3) hold for $G = \overset{n}{\square} K_k$.

This concludes the proof of Theorem A. □

**Proof of Corollary 2**

Let $T \in \mathcal{T}(n)$. Then $T - r$ has $n+1$ connected components, otherwise $T$ contains a cycle, a contradiction. It follows from Lemma 8 that each component of $T - r$ has dimension at most $n - 1$ and thus $T \in \mathcal{H}(n)$. Hence Theorem A implies that there exists $k \in \mathbb{N}$ so that $\overset{n}{\square} K_k \hookrightarrow\!\!\!\!\rightarrow T$. □

## 4 Proof of Theorem B

Let $H \in \mathcal{F}(n)$, i.e., $H = F \cup T$, where $F$ is a graph with $\dim(F) \le n - 1$, $T$ is an $(n+1)$-legged spider of height $n + 2$ such that $V' := V(F) \cap V(T)$ is the set of leaves of $T$. We denote the head of $T$ by $r$ and legs of $T$ by $L_1, L_2, \ldots, L_{n+1}$ and for each $i \in [n+1]$ we let $L_i = p_0(i), p_1(i), \ldots, p_{n+2}(i)$ with $p_0(i) = r$ and $p_{n+2}(i) \in V(F)$. We shall write $p_j = p_j(i)$ when the leg in which $p_j$ is contained is clear from the context.

Let $G = \overset{n}{\square} K_k$ be the Hamming graph on a vertex set $\{0, 1, \ldots, k-1\}^n$, where $k \in \mathbb{N}$ is to be specified later. We shall show that $G \hookrightarrow\!\!\!\!\!\rightarrow H$. By Lemma 10, $G \hookrightarrow\!\!\!\!\!\rightarrow H$ if and only if

S1) $H \not\subseteq_I G$,

S2) for any edge $e \in E(\overline{G})$, $H \subseteq_I G + e$, and

S3) there is an edge $f \in E(\overline{H})$ such that $H + f \subseteq_I G$.

The statement S1) holds since $\{r\} \cup N(r)$ induces a star $K_{1,n+1}$ and therefore by Lemma 8, $H \not\subseteq_I G$.

Now we prove the statement S2). Let $e \in E(\overline{G})$ and $G^+ = G + e$. We need to show that $H \subseteq_I G^+$. Denote the endpoints of $e$ by $\mathbf{u}$ and $\mathbf{v}$, and let $q = \text{dist}_G(\mathbf{u}, \mathbf{v})$, $2 \leq q \leq n$. By Lemma 9, we can suppose that $\mathbf{u} = \mathbf{0}$ and $\mathbf{v}$ has the form $\mathbf{v} = (\underbrace{2, \ldots, 2}_{q}, 0, \ldots, 0)$.

We shall construct an embedding $\sigma$ of $H = F \cup T$ in $G^+$, by first defining $\sigma$ on $V(F)$, then on $r$, and then on the legs of $T$. Let $\sigma$ map $V(F)$ into $\{n+4, n+5, \ldots, k-1\}^n$ such that $[\sigma(v)]_n = k - 1$ for all $v \in V(F)$ and such that $\sigma|_{V(F)}$ is an induced embedding of $F$. Such an embedding exists for sufficiently large $k$, since $\dim(F) \leq n - 1$. Let $\sigma(r) = \mathbf{0}$. Now, we define the embedding $\sigma$ on the legs $L_1, \ldots, L_{n+1}$. First, we give an embedding of $L_1, \ldots, L_{n-1}$, then of $L_n$, and finally of $L_{n+1}$.

Let $i \in [n-1]$ and let $p_j = p_j(i)$ for all $j \in [n+2]$, i.e., $L_i = p_0, p_1, \ldots, p_{n+2}$. Recall that $p_0 = r$, $p_{n+2} \in V(F)$, thus $\sigma(p_{n+2}) = (y_1, y_2, \ldots, y_{n-1}, k-1)$, where $y_j \in \{n+4, \ldots, k-1\}$. We define $\sigma$ on $V(L_i)$ as

$$\sigma(p_1) = (0, 0, 0, \ldots, \ 0, \ \ 1, \ 0, \ \ 0, \ 0, \ldots, \ 0, \ \ \ \ 0 \ ) = \mathbf{e}_i,$$
$$\sigma(p_2) = (y_1, 0, 0, \ldots, \ 0, \ \ 1, \ 0, \ \ 0, \ 0, \ldots, \ 0, \ \ \ \ 0 \ ),$$
$$\sigma(p_3) = (y_1, y_2, 0, \ldots, \ 0, \ \ 1, \ 0, \ \ 0, \ 0, \ldots, \ 0, \ \ \ \ 0 \ ),$$
$$\vdots$$
$$\sigma(p_i) = (y_1, y_2, y_3, \ldots, y_{i-1}, \ 1, \ 0, \ \ 0, \ 0, \ldots, \ 0, \ \ \ \ 0 \ ),$$
$$\sigma(p_{i+1}) = (y_1, y_2, y_3, \ldots, y_{i-1}, \ 1, \ y_{i+1}, \ 0, \ 0, \ldots, \ 0, \ \ \ \ 0 \ ),$$
$$\sigma(p_{i+2}) = (y_1, y_2, y_3, \ldots, y_{i-1}, \ 1, \ y_{i+1}, \ y_{i+2}, \ 0, \ldots, \ 0, \ \ \ \ 0 \ ),$$
$$\vdots$$
$$\sigma(p_{n-1}) = (y_1, y_2, y_3, \ldots, y_{i-1}, \ 1, \ y_{i+1}, \ y_{i+2}, \ldots, \ y_{n-1}, \ \ 0 \ ),$$
$$\sigma(p_n) = (y_1, y_2, y_3, \ldots, y_{i-1}, \ 1, \ y_{i+1}, \ y_{i+2}, \ldots, \ y_{n-1}, \ i+2),$$
$$\sigma(p_{n+1}) = (y_1, y_2, y_3, \ldots, y_{i-1}, \ y_i, \ y_{i+1}, \ y_{i+2}, \ldots, \ y_{n-1}, \ i+2),$$
$$\sigma(p_{n+2}) = (y_1, y_2, y_3, \ldots, y_{i-1}, \ y_i, \ y_{i+1}, \ y_{i+2}, \ldots, \ y_{n-1}, \ k-1).$$

Note that $2 < i + 2 < n + 4 \leq k - 1$ for $i \in [n-1]$.

For $i = n$, let $p_j = p_j(n)$ for all $j \in [n+2]$, i.e., $L_n = p_0, p_1, \ldots, p_{n+2}$ and let $\sigma(p_{n+2}) = (y_1, y_2, \ldots, y_{n-1}, k-1)$, where $y_j \in \{n+4, \ldots, k-1\}$. We define $\sigma$ on $V(L_n)$ as

$$\sigma(p_1) = (0, \ 0, \ \ldots, \ 0, \ \ 0, \ \ 1 \ ) = \mathbf{e}_n,$$
$$\sigma(p_2) = (y_1, \ 0, \ \ldots, \ 0, \ \ 0, \ \ 1 \ ),$$

$$\sigma(p_3) = (y_1, y_2, 0, \ldots, \quad 0, \quad 0, \quad 1 \ ),$$
$$\vdots$$
$$\sigma(p_{n-1}) = (y_1, y_2, \ldots, y_{n-2}, \quad 0, \quad 1 \ ),$$
$$\sigma(p_n) = (y_1, y_2, \ldots, y_{n-2}, \quad 0, \quad n+2),$$
$$\sigma(p_{n+1}) = (y_1, y_2, \ldots, y_{n-2}, y_{n-1}, n+2),$$
$$\sigma(p_{n+2}) = (y_1, y_2, \ldots, y_{n-2}, y_{n-1}, k-1).$$

For $i = n+1$, let $p_j = p_j(n+1)$ for all $j \in [n+2]$, i.e., $L_{n+1} = p_0, p_1, \ldots, p_{n+2}$ and let $\sigma(p_{n+2}) = (y_1, y_2, \ldots, y_{n-1}, k-1)$, where $y_j \in \{n+4, \ldots, k-1\}$. We define $\sigma$ on $V(L_{n+1})$ as

$$\sigma(p_1) = (2, 2, \ldots, 2, 0 \ldots, \quad 0, \quad 0 \ ) = \mathbf{v},$$
$$\sigma(p_2) = (2, 2, \ldots, 2, 0, \ldots, \quad 0, \quad n+3),$$
$$\sigma(p_3) = (2, 2, \ldots, 2, 0, \ldots, y_{n-1}, n+3),$$
$$\vdots$$
$$\sigma(p_n) = (2, y_2, \ldots, y_{n-2}, y_{n-1}, n+3),$$
$$\sigma(p_{n+1}) = (y_1, y_2, \ldots, y_{n-2}, y_{n-1}, n+3),$$
$$\sigma(p_{n+2}) = (y_1, y_2, \ldots, y_{n-2}, y_{n-1}, k-1).$$

**Claim.** $\sigma$ is an induced embedding of $H$ in $G^+$

It is clear from the way we defined $\sigma$ that for any distinct $u, v \in V(H)$, $\sigma(u) \neq \sigma(v)$ and if $uv \in E(H)$, then $\sigma(u)\sigma(v) \in E(G^+)$.

Thus we only need to verify that nonadjacent pairs of vertices of $H$ are mapped to nonadjacent pairs of vertices of $G^+$. Let $x, y \in V(H)$ such that $xy \notin E(H)$ and denote $\sigma(x)$ by $\mathbf{x}$ and $\sigma(y)$ by $\mathbf{y}$. To prove that $\mathbf{xy} \notin E(G^+)$, it suffices to show that $\mathbf{x}$ and $\mathbf{y}$ differ in at least two coordinates.

Let $x, y \in V(F)$, then $\mathbf{xy} \notin E(G^+)$ since $\sigma|_{V(F)}$ is an induced embedding of $F$.

Let $x \notin V(F)$, $y \in V(F)$ and $x$ has a neighbor $x'$ in $V(F)$. By construction, $\mathbf{x}$ differs from $\mathbf{x'} := \sigma(x')$ in exactly the last coordinate. Moreover, $\mathbf{y}$ differs from $\mathbf{x'}$ in at least one coordinate, denote it by $j$. Note that $j \neq n$ since $[\mathbf{x'}]_n = [\mathbf{y}]_n = k - 1$. Thus $\mathbf{x}$ differs from $\mathbf{y}$ in the $j^{\text{th}}$ and the $n^{\text{th}}$ coordinates.

Let $x \notin V(F)$, $y \in V(F)$ and $x$ does not have any neighbors in $V(F)$. Let $i \in [n+1]$ such that $x \in V(L_i)$. If $i \in [n-1]$, then $[\mathbf{x}]_i = 1 \neq [\mathbf{y}]_i$, if $i = n$, then $[\mathbf{x}]_{n-1} = 0 \neq [\mathbf{y}]_{n-1}$, and if $i = n+1$, then $[\mathbf{x}]_1 = 2 \neq [\mathbf{y}]_1$. Moreover for any value of $i$, $[\mathbf{x}]_n \in \{0, 1, i+2\} \not\ni [\mathbf{y}]_n$, thus $\mathbf{x}$ and $\mathbf{y}$ differ in at least two coordinates.

Let $x, y \notin V(F)$ and $x, y$ are in the same leg. It is clear that $\mathbf{xy} \notin E(G^+)$.

Let $x \in V(L_i) \setminus V(F)$ and $y \in V(L_j) \setminus V(F)$ with $1 \leq i < j \leq n+1$, $x \neq r$, $y \neq r$. If one

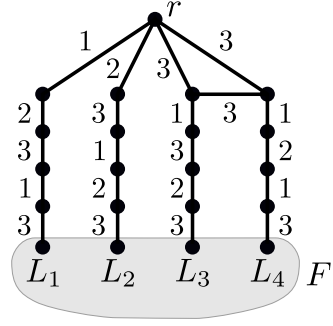

Figure 4: The coloring $c$ of the legs $L_1, L_2, L_3, L_4$ when $H \in \mathcal{F}(3)$.

of them, say $x$, has a neighbor in $V(F)$, i.e., $x = p_{n+1}(i)$, then $[\mathbf{x}]_n = i + 2 \neq [\mathbf{y}]_n$. Thus if $\mathbf{xy} \in E(G^+)$, then $\mathbf{x}$ and $\mathbf{y}$ agree on the first $n-1$ coordinates, from which $p_{n+2}(i) = p_{n+2}(j)$ follows, a contradiction. Hence $\mathbf{xy} \notin E(G^+)$. Finally, assume that neither $x$ nor $y$ has a neighbor in $V(F)$. If $j \leq n$, then $[\mathbf{x}]_i = 1 \neq [\mathbf{y}]_i$ and $[\mathbf{y}]_j \in \{1, n+2\} \not\ni [\mathbf{x}]_j$. If $j = n+1$, then $[\mathbf{y}]_1 = 2 \neq [\mathbf{x}]_1$ and either $[\mathbf{y}]_2 = 2 \neq [\mathbf{x}]_2$ or $[\mathbf{y}]_n = n+3 \neq [\mathbf{x}]_n$.

This proves the claim and the statement S2).

To prove statement S3), we let $f = xy$, where $x = p_1(n)$ and $y = p_1(n+1)$, i.e., $x$ and $y$ are neighbors of $r$ on legs $L_n$ and $L_{n+1}$ respectively. We shall show that $G$ contains an induced copy of $H^+ = H + f$, using the coloring criterium (Theorem 7).

Since $\dim(F) \leq n-1$, there is a nice coloring $c'$ of $E(F)$ with colors from $[n-1]$. We shall define a coloring $c$ of $E(H^+)$ with colors from $[n]$. Let $c(e) = c'(e)$ for all $e \in E(F)$.

Let $c(L_i)$ be the list of colors on the leg $L_i$ starting with an edge incident to $r$, i.e., $c(L_i) = \bigl(c(e_1(i)), c(e_2(i)), \ldots, c(e_{n+2}(i))\bigr)$ where $e_j(i) = p_{j-1}(i)p_j(i)$ for $j \in [n+2]$ are the edges of $L_i$. For $n = 3$, we illustrate the coloring of the legs $L_1, L_2, L_3, L_4$ of $T$ in Figure 4. For $n \geq 4$, define the coloring $c$ on the legs as follows:

$$
\begin{aligned}
c(L_1) &= (\ 1, \quad 2, \quad \ldots, n-1, \quad n, \quad 1, \quad n), \\
c(L_2) &= (\ 2, \quad 3, \quad \ldots, \quad n, \quad 1, \quad 2, \quad n), \\
&\quad \vdots \\
c(L_{n-1}) &= (n-1, n, 1, \ldots, n-3, n-2, n-1, n), \\
c(L_n) &= (\ n, \quad 1, 2, \ldots, n-2, \quad n, \quad n-1, n), \\
c(L_{n+1}) &= (\ n, \quad 1, 2, \ldots, n-2, n-1, \quad 1, \quad n).
\end{aligned}
$$

Color $xy$ with $n$, so, we have $c(rx) = c(ry) = c(xy) = n$. Note that each leg contains all the colors from $[n]$ and ends with the same color $n$.

**Claim.** The coloring $c$ is a nice coloring.

The vertex set $\{r, x, y\}$ induces a triangle of color $n$. Any other triangle is in $F$ and is monochromatic as $c'$ is a nice coloring.

Let $u, v$ be any two distinct, non-adjacent vertices of $H^+$. We need to show that there are two distinct colors that appear on each induced $u, v$-path.

Case 1. $u, v \in V(F)$.

Since $c'$ is a nice coloring, there is a pair of distinct colors $i, j$ appearing on all induced $u, v$-paths in $F$. Any other induced $u, v$-path in $H^+$ contains all the edges of two legs of $T$ (except for maybe $rx, ry$) and so it contains all the colors from $[n]$. Thus the colors $i$ and $j$ appear on such a path as well.

Case 2. $u, v \in V(L_j)$, $j \in [n+1]$.

Any induced $u, v$-path is either a subpath $P$ of $L_j$ or contains all the edges of a leg $L_i$, $i \neq j$ (except for maybe $rx, ry$) and thus it contains all the colors from $[n]$. Since $P$ contains at least two distinct colors, these two colors are present on each induced $u, v$-path.

Case 3. $u, v \notin V(F)$, $u \in V(L_i)$, $v \in V(L_j) \setminus V(L_i)$, $1 \leq i < j \leq n+1$.

Then $v \neq r$ and we can assume that $u \neq r$, otherwise we are in Case 2. Let $u' \in V(L_i) \cap V(F)$, $v' \in V(L_j) \cap V(F)$ be leaves of $L_i$ and $L_j$, respectively. Let $P'_u$ be the subpath of $L_i$ from $u$ to $r$ and $P''_u$ be subpath of $L_i$ from $u$ to $u'$. Similarly, let $P'_v$ be the subpath of $L_j$ from $v$ to $r$ and $P''_v$, subpath of $L_j$ from $v$ to $v'$. See Figure 5.

If an induced $u, v$-path contains an edge of a leg that is not in $E(L_i \cup L_j)$, then it is easy to see that this path contains all the colors from $[n]$. So, we can consider only $u, v$-paths that are either equal to $P'_u \cup P'_v$ (or $((P'_u \cup P'_v) - r) \cup \{xy\}$ if $i = n$, $j = n+1$), or a union of $P''_u \cup P''_v$ and a $u', v'$-path in $F$.

Assume first that $u$ is adjacent to $u'$ or $v$ is adjacent to $v'$.

Without loss of generality, $u$ is adjacent to $u'$. Then $P'_u$ contains all the colors from $[n]$ and $P''_u$ contains the color $n$. By the niceness of $c'$ and Observation 6, there is a color $\ell \in [n-1]$ that appears on each induced $u', v'$-path in $F$ and so the union of $P''_u \cup P''_v$ and any $u', v'$-path in $F$ contains the colors $\ell, n$. Thus all induced $u, v$-paths contain the colors $\ell$ and $n$.

Assume that $u$ is not adjacent to $u'$ and $v$ is not adjacent to $v'$.

If $j < n+1$, then $P'_u, P''_u$ contain color $i$ and $P'_v, P''_v$ contain color $j$, thus all induced $u, v$-paths contain the colors $i$ and $j$.

If $j = n+1$ and $v$ is not adjacent to $r$, then each of $P'_v$ and $P''_v$ contains the colors $1, n$ and so the colors $1, n$ appear on any induced $u, v$-path.

If $j = n+1$ and $v$ is adjacent to $r$, we distinguish two cases based on $i$. If $i < n$, then each of $P'_u \cup P'_v$ and $P''_u \cup P''_v$ contains the colors $i, n$ and thus all induced $u, v$-paths contain the colors $i$ and $j$. If $i = n$, then $u$ is not adjacent to $r$, otherwise $uv \in E(H^+)$, a contradiction. Hence the colors $1, n$ appear on each of $P'_u \cup P'_v$ and $P''_u \cup P''_v$ and so all induced $u, v$-paths contain the colors $1$ and $n$.

Case 4. $v \in V(L_j)$, $j \in [n+1]$ and $u \in V(F)$.

We can assume that $v \notin V(F)$ and $u \notin V(L_j)$ otherwise we are in Case 1 or in Case 2, respectively. Let $v' \in V(L_j) \cap V(F)$ be a leaf of $L_j$ and $P$ be the subpath of $L_j$ from $v$ to $v'$. Then $P$ contains the color $n$. By the niceness of $c'$ and Observation 6, there is a color $\ell \in [n-1]$ that appears on each induced $u, v'$-path in $F$. Let $P'$ be an induced $u, v$-path that is not the union of $P$ and an induced $u, v'$-path in $F$. Then $P$ contains all the edges of a leg $L_i$, $i \neq j$ (except for maybe $rx, ry$) and thus it contains all the colors from $[n]$. Hence the colors $\ell$ and $n$ appear on all induced $u, v$-paths.

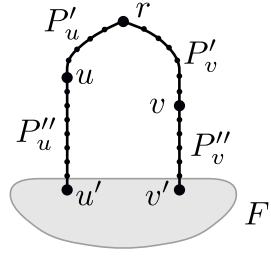

Figure 5: Illustration for the proof of Theorem B.

This concludes the proof of Claim, Statement S3), and the proof of Theorem B. $\square$

**Proof of Corollary 3.** Let $v \in V(H)$ be the (unique) vertex of maximum degree $n+1$, $F = H - v$ and $F' \subseteq H'$ be an $(n+1)$-subdivision of $F$, so $V(F) \subseteq V(F')$. To apply Theorem B, we need to verify that $H' \in \mathcal{F}(n)$. For that it suffices to show that $F'$ has dimension at most $n-1$. We prove this by defining a nice coloring $c : E(F') \to [n-1]$ and using the coloring criterium (Theorem 7). The maximum degree of $F'$ is at most $n-1$. Denote the set of paths of length $n+2$ that have both endpoints in $V(F)$ by $\mathcal{P}$, i.e., $\mathcal{P}$ is the set of paths that arise from the subdivision of an edge of $F$. Color all edges of $F'$ incident to vertices corresponding to vertices of $F$ with at most $n-1$ colors properly. This gives a coloring of the end-edges of each path in $\mathcal{P}$. For a path $P = v_1, v_2, \ldots, v_{n+2}$ and an edge-coloring $c$, we use the notation $c(P) = \big(c(v_1v_2), c(v_2v_3), \ldots, c(v_{n+1}v_{n+2})\big)$. Let $P \in \mathcal{P}$. If the end-edges of $P$ have different colors $i, j$ with $i > j$, then let $c(P) = (i, 1, 2, 3, \ldots, n-1, x, j)$, where $x \neq j$ and $x \neq n-1$ (such an $x$ exists, since $n-1 \geq 3$). If the end-edges of $P$ have the same color $i$, then let $c(P) = (i, 1, 2, 3, \ldots, n-1, 1, i)$ if $i \neq 1$ and $c(P) = (1, 2, 1, 3, 4, \ldots, n-1, 2, 1)$ if $i = 1$. Note that any path in $\mathcal{P}$ contains all the colors from $[n-1]$.

Now we verify that $c$ is a nice coloring. Since $F'$ is triangle-free, we only need to check that for any two nonadjacent vertices $u, v$, there are distinct colors $i, j$ that appear in any induced $u, v$-path. Let $P_u, P_v$ be some paths in $\mathcal{P}$ such that $u \in V(P_u)$ and $v \in V(P_v)$.

If $V(P_u) \cap V(P_v) = \emptyset$, then any $u, v$-path contains an entire path from $\mathcal{P}$ and thus all the colors from $[n-1]$.

If $V(P_u) \cap V(P_v) = \{w\}$, then $w$ is an endpoint of both paths, that is, it corresponds to a vertex of $F$ and so the end-edges of $P_u$ and $P_v$ that are incident to $w$ are colored with distinct colors, $i, j$. Any other $u, v$-path contains an entire path from $\mathcal{P}$ and thus all the colors from $[n-1]$.

If $|V(P_u) \cap V(P_u)| \geq 2$, then $P_u = P_v$ and since $c$ is a proper coloring, there are distinct colors $i, j$ that appear in the shortest $u, v$-path (which is a subpath of $P_u = P_v$). Any other $u, v$-path contains an entire path from $\mathcal{P}$ and thus all the colors from $[n-1]$. $\square$

# 5 Graphs induced-saturated in the Cartesian product of two cliques.

The aim of this section is to give a complete characterization of connected graphs $H$ so that $K_k \square K_k \hookrightarrow\!\!\!\!\!/\, H$ for some $k \in \mathbb{N}$:

**Theorem 12.** *Let $H$ be a connected graph. There is $k \in \mathbb{N}$ such that $K_k \square K_k \hookrightarrow\!\!\!\!\!/\, H$ if and only if $H$ is either chipped or $H \in \mathcal{A} \cap \mathcal{B}$.*

To define the notion of chipped graphs and classes $\mathcal{A}$ and $\mathcal{B}$ we need some definitions.

## 5.1 Definitions

We use the notation $|H|$ for the number of vertices of a graph $H$ and $\{\{X\}\}$ for the complete graph on the vertex set $X$. We call a graph $H$ *2-Hamming* if there exists $k \in \mathbb{N}$ so that $H \subseteq_I K_k \square K_k$. We say that $\{F_1, F_2\}$ is a *2-Hamming decomposition* of a graph $H$ if $V(F_1) = V(F_2) = V(H)$, $\{E(F_1), E(F_2)\}$ is a partition of $E(H)$, and $F_i$ is a pairwise vertex-disjoint union of cliques for $i \in \{1, 2\}$. Note that a clique in $F_1$ intersects any clique in $F_2$ in at most one vertex. The coloring criterium (Theorem 7) implies that a graph is 2-Hamming if and only if it has a nice coloring with two colors. It is easy to see that a 2-Hamming decomposition corresponds to a nice coloring with colors 1 and 2, such that the edges of $F_1$ have color 1 and the edges of $F_2$ have color 2. Thus a graph is 2-Hamming if and only if it has a 2-Hamming decomposition. One can show that $\min\{k : H \subseteq_I K_k \square K_k\}$ is equal to the smallest integer $\ell$ such that $H$ has a 2-Hamming decomposition $\{F_1, F_2\}$ such that each of $F_1$ and $F_2$ has at most $\ell$ connected components. We prove this fact formally in Lemma 17. For a graph $G$, let $\mathcal{C}(G)$ denote the number of connected components of $G$.

**Definition 13.** *Let $H$ be a 2-Hamming graph. We say that a 2-Hamming decomposition $\{F_1, F_2\}$ of $H$ is **optimal**, if $\min\{k : H \subseteq_I K_k \square K_k\} = \max_{i \in \{1,2\}} \mathcal{C}(F_i)$. We call $H$ **balanced** if $\mathcal{C}(F_1) = \mathcal{C}(F_2)$ for an optimal 2-Hamming decomposition $\{F_1, F_2\}$ of $H$ and **unbalanced** otherwise.*

**Definition of chipped graphs**

We call a graph $H$ *chipped* if $H = H' - e$, where $H'$ is a connected, 2-Hamming graph with a 2-Hamming decomposition $\{F_1, F_2\}$, such that $F_1$ has a component $K$ of order 3 or 4 that contains the endpoints $u, v$ of $e$ and one of the following holds:

ch1) $V(K) = \{u, v, w\}$, $\deg_{H'}(w) \geq 3$;

ch2) $V(K) = \{u, v, w\}$ and $H' - E(K)$ has two components, one of which is $\{\{w\}\}$;

ch3) $V(K) = \{u, v, w, x\}$, $\deg_{H'}(w) = 3$ and $H' - E(K)$ has a connected component $\widetilde{H}$ so that $V(K) \cap V(\widetilde{H}) = \{u\}$.

See Figure 6 for an illustration.

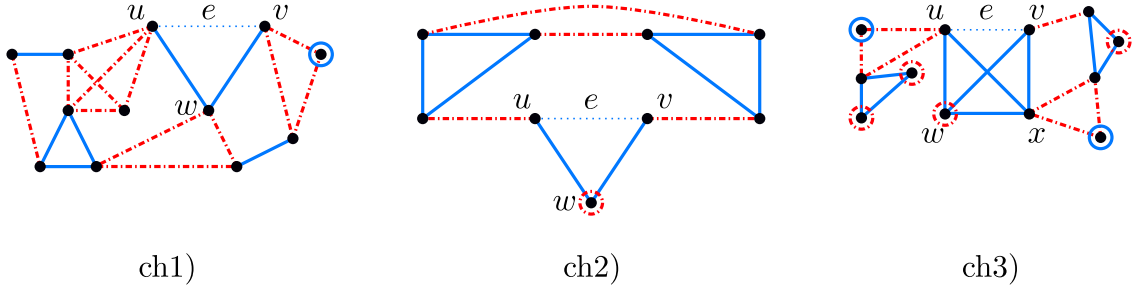

ch1)  ch2)  ch3)

Figure 6: Examples of chipped graphs. Colors indicate a nice 2-coloring of $H'$, each color also corresponds to edges and singleton cliques of $F_i$, $i \in \{1, 2\}$.

**Definition of the graph classes $\mathcal{A}$ and $\mathcal{B}$**

Let $H$ be a 2-Hamming graph with 2-Hamming decomposition $\{F_1, F_2\}$. For a subgraph $H'$ of $H$ and $i \in \{1, 2\}$, let $F_i|H'$ denote the subgraph of $F_i$ that is the union of those components of $F_i$ whose vertex sets are completely contained in $V(H')$.

Let $\mathcal{A}$ be the class of connected, unbalanced, 2-Hamming graphs $H$ that have an optimal 2-Hamming decomposition $\{F_1, F_2\}$ with $\mathcal{C}(F_1) > \mathcal{C}(F_2)$ such that $H$ has two nonadjacent vertices $u, v$, for which one of the following holds:

a1) $\deg_{F_1}(u) = \deg_{F_1}(v) = 0$;

a2) there is a vertex $w \in N(u) \cap N(v)$ such that $H - \{uw, wv\}$ has 3 components $H_1, H_2, H_3$ with $\mathcal{C}(F_1|H_1) > \mathcal{C}(F_2|H_1)$, $\mathcal{C}(F_1|H_2) > \mathcal{C}(F_2|H_2)$, and $H_3 = \{\{w\}\}$.

See Figure 7 for an illustration.

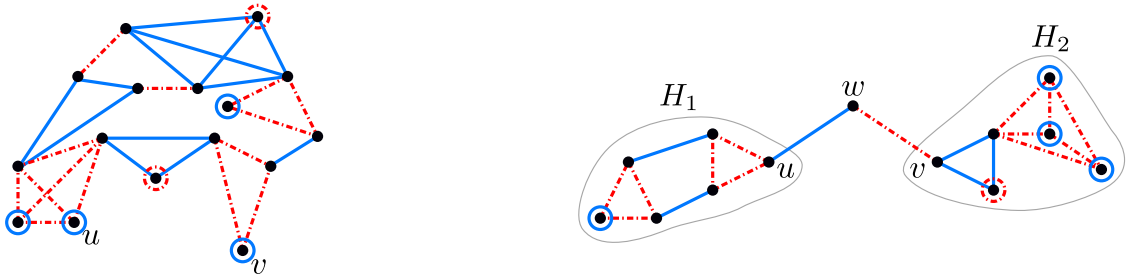

a1) $\mathcal{C}(F_1) = 8$ and $\mathcal{C}(F_2) = 7$.    a2) $\mathcal{C}(F_1) = 8$, $\mathcal{C}(F_2) = 5$, $\mathcal{C}(F_1|H_1) = 3$, $\mathcal{C}(F_2|H_1) = 2$, $\mathcal{C}(F_1|H_2) = 4$, $\mathcal{C}(F_2|H_2) = 2$.

Figure 7: Examples of graphs from class $\mathcal{A}$.

Let $\mathcal{B}$ be the class of connected, unbalanced, 2-Hamming graphs $H$ that have an optimal 2-Hamming decomposition $\{F_1, F_2\}$ with $\mathcal{C}(F_1) > \mathcal{C}(F_2)$ such that $H$ has two adjacent vertices $u, v$ for which one of the following holds:

b1) the edge $e = uv$ is a bridge with $e \in E(F_i)$ so that the components $H_1$, $H_2$ of $H - e$ satisfy $\mathcal{C}(F_1|H_1) \geq \mathcal{C}(F_2|H_1) + i$ and $\mathcal{C}(F_1|H_2) \geq \mathcal{C}(F_2|H_2) + i$.

b2) there is a vertex $w \in N(u) \cap N(v)$, i.e., $\{u, v, w\}$ induce a triangle $T$ in $H$, and $H - E(T)$ has 3 components $H_1$, $H_2$, $H_3$ with $\mathcal{C}(F_1|H_1) > \mathcal{C}(F_2|H_1)$, $\mathcal{C}(F_1|H_2) > \mathcal{C}(F_2|H_2)$, and $H_3 = \{\{w\}\}$.

See Figure 8 for an illustration.

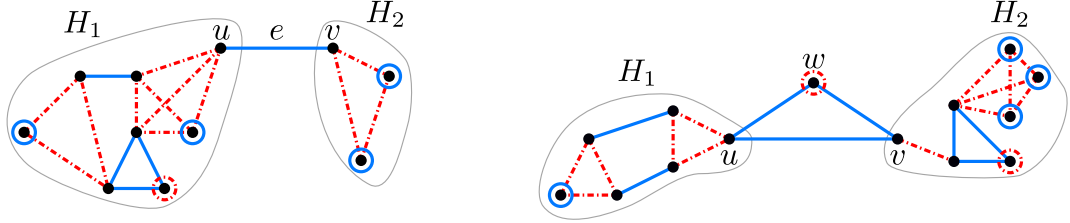

b1) $\mathcal{C}(F_1) = 7$, $\mathcal{C}(F_2) = 4$, $\mathcal{C}(F_1|H_1) = 4$, $\mathcal{C}(F_2|H_1) = 3$, $\mathcal{C}(F_1|H_2) = 2$, $\mathcal{C}(F_2|H_2) = 1$.

b2) $\mathcal{C}(F_1) = 8$, $\mathcal{C}(F_2) = 6$, $\mathcal{C}(F_1|H_1) = 3$, $\mathcal{C}(F_2|H_1) = 2$, $\mathcal{C}(F_1|H_2) = 4$, $\mathcal{C}(F_2|H_2) = 3$.

Figure 8: Examples of graphs from class $\mathcal{B}$.

## 5.2 Preliminary Lemmas

**Lemma 14.** *Let $k \in \mathbb{N}$, $G = K_k \square K_k$ and $H$ be a graph. Then $G \hookrightarrow\!\!\!\!\!\twoheadrightarrow H$ if and only if*

s1) $H \not\subseteq_I G$,

s2) *there is an edge $e \in E(H)$ such that $H - e \subseteq_I G$ and*

s3) *there is an edge $f \in E(\overline{H})$ such that $H + f \subseteq_I G$.*

*Proof.* By Lemma 10, $G \hookrightarrow\!\!\!\!\!\twoheadrightarrow H$ if and only if
S1) $H \not\subseteq_I G$,
S2) for any $e \in E(\overline{G})$, $H \subseteq_I G + e$, and
S3) there is an edge $f \in E(\overline{H})$ such that $H + f \subseteq_I G$.
Thus it is sufficient to prove that if $H \not\subseteq_I G$, then s2) is equivalent to S2).

Assume that $H \not\subseteq_I G$ and S2) holds. Let $e = uv \in E(\overline{G})$ and let $\sigma : V(H) \to V(G)$ be an induced embedding of $H$ into $G + e$. Since $H \not\subseteq_I G$, the vertex set $\sigma(V(H))$ has to contain both $u$ and $v$. Since $\sigma$ is an induced embedding, $\sigma^{-1}(u)\sigma^{-1}(v) \in E(H)$. Let $e' = \sigma^{-1}(u)\sigma^{-1}(v)$, then $\sigma$ is an induced embedding of $H - e'$ into $G$, i.e., s2) holds.

Now assume that $H \not\subseteq_I G$ and s2) holds. Let $e' = u'v' \in E(H)$ so that $H - e' \subseteq_I G$ and let $\sigma : V(H) \to V(G)$ be an induced embedding of $H - e'$ into $G$. Since $\sigma$ is an induced embedding, $\sigma(u')\sigma(v') \notin E(G)$. Let $f = \sigma(u')\sigma(v')$, then $\sigma$ is an induced embedding of $H$ into $G + f$, i.e., there exists an edge $f \in E(\overline{G})$ such that $H \subseteq_I G + f$. Since the distance of two distinct vertices in $G$ is either 1 or 2, any pair of nonadjacent vertices have distance 2. Thus Lemma 9 implies that if there is an edge $f \in E(\overline{G})$ such that $H \subseteq_I G + f$, then $H \subseteq_I G + e$ for any $e \in E(\overline{G})$, i.e., S2) holds. □

**Lemma 15.** *Let $H$ be a 2-Hamming graph, $\{F_1, F_2\}$ be a 2-Hamming decomposition of $H$ and $K$ be a component in $F_1$. Then for any $v \in V(H) \setminus V(K)$, $v$ is adjacent to at most one vertex of $K$. In particular, for any vertex/edge of $H$, any maximal clique in $H$ containing this vertex/edge is a component of either $F_1$ or $F_2$.*

*Proof.* Suppose that $v \in V(H) \setminus V(K)$ is adjacent to two distinct vertices $u, u' \in V(K)$. Then $uu' \in E(F_1)$ and $uv, u'v \in E(F_2)$. Thus $u$ and $u'$ are in the same connected component of $F_2$. Since each component of $F_2$ is a clique, $uu' \in E(F_2)$, a contradiction since $uu' \in E(F_1)$. □

**Lemma 16.** *Let $k \geq 2$. A graph that contains an induced copy of $K_{1,3}$, $K_4 - e$, or $C_{2k+1}$ is not 2-Hamming.*

*Proof.* It is easy to see that the listed graphs do not have a nice coloring with 2 colors, thus by the coloring criterium (Theorem 7) their dimension is at least 3. Therefore, any graph that contains an induced copy of them has dimension at least 3, i.e., is not 2-Hamming. □

**Lemma 17.** *For any 2-Hamming graph $H$,*

$$\min\{k : H \subseteq_I K_k \square K_k\} = \min_{\{F_1, F_2\}} \max_{i \in \{1,2\}} \mathcal{C}(F_i),$$

*where the minimum on the right hand side is taken over all $\{F_1, F_2\}$ 2-Hamming decompositions of $H$.*

*Proof.* Let $c = \min \{ \max\{\mathcal{C}(F_1), \mathcal{C}(F_2)\} : \{F_1, F_2\} \text{ is a 2-Hamming decomposition of } H\}$ and $k$ be the smallest integer such that $H \subseteq_I K_k \square K_k$.

We show that $k \leq c$ by proving that $H \subseteq_I K_c \square K_c$. Let $\{F_1, F_2\}$ be a 2-Hamming decomposition of $G$ so that $c = \max_{i \in \{1,2\}} \mathcal{C}(F_i)$. Without loss of generality, $c = \mathcal{C}(F_1) \geq \mathcal{C}(F_2)$. Let $C_1^1, C_2^1, \ldots, C_c^1$ be the components of $F_1$ and $C_1^2, C_2^2, \ldots, C_d^2$ be the components of $F_2$ with $d = \mathcal{C}(F_2) \leq c$. Since $F_1$ and $F_2$ are spanning subgraphs, each vertex of $G$ belongs to exactly one component of $F_1$ and exactly one component of $F_2$.

We define an embedding $\sigma$ of $H$ into $K_c \square K_c$ as follows: for any $v \in V$, let $\sigma(v) = (i, j)$ if $v \in C_i^1$ and $v \in C_j^2$. Since any component of $F_1$ intersects a component of $F_2$ in at most one vertex, $\sigma$ is injective. Moreover, $\sigma$ is an induced embedding, since

$$uv \in E(H) \iff uv \in E(F_1) \text{ or } uv \in E(F_2)$$
$$\iff u \text{ and } v \text{ are in the same component either in } F_1 \text{ or in } F_2$$
$$\iff \text{dist}_{\text{Ham}}(\sigma(u), \sigma(v)) = 1 \iff \sigma(u)\sigma(v) \in E(K_c \square K_c).$$

We show that $k \geq c$ by defining a 2-Hamming decomposition $\{F_1, F_2\}$ of $H = (V, E)$ such that $\max\{\mathcal{C}(F_1), \mathcal{C}(F_2)\} \leq k$. Let $\sigma$ be an induced embedding of $H$ into $K_k \square K_k$. Let $E = E_1 \cup E_2$ such that $uv \in E_i$ if $\sigma(u)$ and $\sigma(v)$ differ exactly in the $i^{\text{th}}$ coordinate, $i \in \{1, 2\}$. Then $F_1 = (V, E_1)$ and $F_2 = (V, E_2)$ form a 2-Hamming decomposition of $H$. We see that two vertices are nonadjacent in $F_i$ if and only if they differ in the $(2-i)^{\text{th}}$ coordinate. Therefore there are at most $k$ pairwise nonadjacent vertices in $F_i$, i.e. the size of the maximal independent set in $F_i$ is at most $k$. It implies that $F_i$ has at most $k$ connected components, thus $\max_{i \in \{1,2\}} \mathcal{C}(F_i) \leq k$ and so $c \leq k$. □

The following lemma is a consequence of Lemma 15 and the coloring criterium (Theorem 7).

**Lemma 18.** *The 2-Hamming decomposition of a connected graph is unique. In particular, if a 2-Hamming graph $H$ has maximal cliques $K$, $K'$ and two 2-Hamming decompositions $\{F_1, F_2\}$, $\{F_1^*, F_2^*\}$ such that $K$ and $K'$ belong to the same family in one and to different families in the other, i.e., $K, K' \subseteq F_i$ and $K \subseteq F_i^*$, $K' \subseteq F_{3-i}^*$ then $H$ is disconnected and $K$ and $K'$ belong to distinct connected components of $H$.*

## 5.3 Proof of Theorem 12

We call a graph $H$ *good* if $K_k \square K_k \hookrightarrow\!\!\!\!\!\to H$ for some $k \in \mathbb{N}$.

Recall that Theorem 12 claims that a connected graph $H$ is good if and only if it is either chipped or contained in $\mathcal{A} \cap \mathcal{B}$.

Let $H$ be a connected graph. We shall prove a more precise statement: if $\dim(H) > 2$, then $H$ is good if and only if $H$ is chipped; if $\dim(H) \leq 2$, then $H$ is good if and only if $H \in \mathcal{A} \cap \mathcal{B}$.

**Claim 1. If $H$ is chipped then it is good.**
Let $H$ be chipped with $H'$, $e$, $K$, $u, v, w, x$ and $\widetilde{H}$ as in the definition of chipped graphs. Let $\{F_1', F_2'\}$ be a 2-Hamming decomposition of $H'$ and assume, without loss of generality, that $uv \in E(F_1')$. To show that $H$ is good, it is sufficient to prove that $H$ is not 2-Hamming and that there are suitable edges $e \in E(H)$, $f \in E(\overline{H})$ so that $\dim(H - e) \leq 2$ and $\dim(H + f) \leq 2$. Indeed, if $\dim(H - e) \leq 2$ and $\dim(H + f) \leq 2$, then there are integers $l, m \in \mathbb{N}$ so that $H - e \subseteq_I K_l \square K_l$, $H + f \subseteq_I K_m \square K_m$, hence for any $k \geq l, m$ we have $H \not\subseteq_I K_k \square K_k$, $H - e \subseteq_I K_k \square K_k$, and $H + f \subseteq_I K_k \square K_k$, thus by Lemma 14, $K_k \square K_k \hookrightarrow\!\!\!\!\!\to H$. First we show that $H$ is not 2-Hamming by finding an induced copy of a not 2-Hamming graph in $H$ and applying Lemma 16.

If ch1) holds, then by Lemma 15, $H$ contains an induced a star $K_{1,3}$ centered at $w$. If ch2) holds, consider a shortest $u, v$-path $P$ in $H' - E(K)$. Then $P$ is an induced path and so it is $F_1'$-$F_2'$-alternating with its first and last edge in $E(F_2')$. Hence $P$ is an induced odd $u, v$-path and thus $H$ contains an induced odd cycle of length at least 5 formed by the edges $vw, wu$ and the edges of $P$. If ch3) holds, then $V(K)$ induces a graph isomorphic to $K_4$ minus an edge. We see by Lemma 16, that in any of the above cases $H$ is not 2-Hamming.

By the definition of chipped graphs, adding the edge $f = uv$ to $H$ results the 2-Hamming graph $H'$.

If ch1) or ch2) holds, then deleting the edge $uw$ results a 2-Hamming graph. Indeed, $F_1^- = (F_1' - K) \cup \{\{v, w\}\} \cup \{\{u\}\}$, $F_2^- = F_2'$ form a 2-Hamming decomposition of $H - uw$. If ch3) holds, then deleting $ux$ results a 2-Hamming graph since $F_1^- = (F_1' - K - (F_1'|\widetilde{H})) \cup (\{\{w, v, x\}\} \cup (F_2'|\widetilde{H}))$, $F_2^- = (F_2' - \{\{w\}\} - (F_2'|\widetilde{H})) \cup (\{\{u, w\}\} \cup (F_1'|\widetilde{H}))$ form a 2-Hamming decomposition of $H - ux$, see Figure 9.

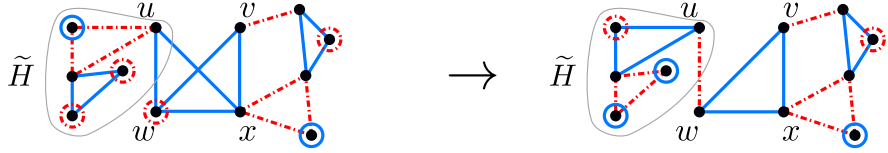

Figure 9: Illustration for the proof of Theorem 12 (Claim 1, case ch3)).

Thus $H$ is good.

**Claim 2. If $H \in \mathcal{A} \cap \mathcal{B}$ then $H$ is good.**

If $H \in \mathcal{A} \cap \mathcal{B}$, then $H$ is an unbalanced 2-Hamming graph. Let $k$ be the smallest integer such that $H \subseteq_I K_{k+1} \square K_{k+1}$.
We shall show that $K_k \square K_k \hookrightarrow H$ by using Lemma 14. In particular, we shall prove that
s1) $H \not\subseteq_I K_k \square K_k$,
s2) there is an edge $e \in E(H)$ such that $H - e \subseteq_I K_k \square K_k$, and
s3) there is an edge $f \in E(\overline{H})$ such that $H + f \subseteq_I K_k \square K_k$.

By the choice of $k$, $H \not\subseteq_I K_k \square K_k$, so s1) follows. We shall prove that s2) is implied by the fact that $H \in \mathcal{B}$ and s3) is implied by the fact that $H \in \mathcal{A}$. Let $\{F_1, F_2\}$ be an optimal 2-Hamming decomposition of $H$ with $\mathcal{C}(F_1) = k + 1$ and $\mathcal{C}(F_2) \leq k$, see Lemma 17.

Let $u, v, H_2$ be as in the definition of the classes $\mathcal{A}$ and $\mathcal{B}$. The rough idea of the proof of this claim is that if we add or delete the edge $uv$, then we can "swap the decomposition in $H_2$", i.e., exchange $F_1|H_2$ and $F_2|H_2$. Then using the unbalancedness assumptions on $\{F_1|H_2, F_2|H_2\}$, we show that the resulting decomposition has at most $k$ components in both families and thus by Lemma 17, $K_k \square K_k$ contains an induced copy of $H + f$ or $H - e$ respectively.

Let $H \in \mathcal{A}$. Let $f = uv$, where $u, v$ are vertices as in the definition of class $\mathcal{A}$ and let $H^+ = H + f$. We shall construct a 2-Hamming decomposition $\{F_1^+, F_2^+\}$ of $H^+$.

If a1) holds, then let $F_2^+ = F_2$ and let $F_1^+$ be obtained from $F_1$ by replacing the singleton components induced by $u$ and $v$ with $\{\{u, v\}\}$. Then $\mathcal{C}(F_1^+) = \mathcal{C}(F_1) - 1 = k$ and $\mathcal{C}(F_2^+) = \mathcal{C}(F_2) \leq k$, implying that $H + f \subseteq_I K_k \square K_k$.

Now assume that a2) holds and let $w, H_1$ and $H_2$ be as in the definition of class $\mathcal{A}$. Then by Lemma 15, $\{\{u,w\}\}$ and $\{\{v,w\}\}$ are a components of $F_1$ or $F_2$. Without loss of generality, $\{\{u,w\}\}$ is a component of $F_1$ and $\{\{v,w\}\}$ is a component of $F_2$. Define $F_1^+ = (F_1|H_1) \cup (F_2|H_2) \cup \{\{u,v,w\}\}$ and $F_2^+ = (F_2|H_1) \cup (F_1|H_2) \cup \{\{w\}\}$, see Figure 10. Then $\mathcal{C}(F_1^+) = \mathcal{C}(F_1|H_1) + \mathcal{C}(F_2|H_2) + 1 < \mathcal{C}(F_1|H_1) + \mathcal{C}(F_1|H_2) + 1 = \mathcal{C}(F_1) = k+1$ and $\mathcal{C}(F_2^+) = \mathcal{C}(F_2|H_1) + \mathcal{C}(F_1|H_2) + 1 < \mathcal{C}(F_1|H_1) + \mathcal{C}(F_1|H_2) + 1 = \mathcal{C}(F_1) = k+1$. Thus $\max\{\mathcal{C}(F_1^+), \mathcal{C}(F_2^+)\} \leq k$, implying that $H + f \subseteq_I K_k \square K_k$.

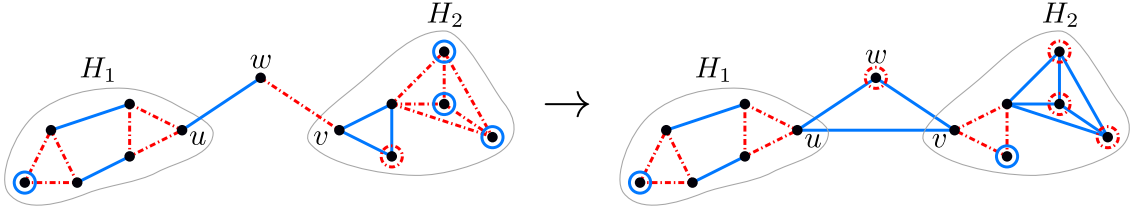

Figure 10: Illustration for the proof of Theorem 12 (Claim 2, case of $H^+$, a2)).

Let $H \in \mathcal{B}$. Let $e = uv$, where $u, v$ are vertices as in the definition of class $\mathcal{B}$ and let $H^- = H - e$. We shall construct a 2-Hamming decomposition $\{F_1^-, F_2^-\}$ of $H^-$.

Let b1) hold for $H$ and let $H_1$, $H_2$ be the components of $H - e$. Then without loss of generality, $u \in V(H_1)$ and $v \in V(H_2)$.
If $e \in E(F_1)$, we have $\mathcal{C}(F_1|H_1) \geq \mathcal{C}(F_2|H_1) + 1$ and $\mathcal{C}(F_1|H_2) \geq \mathcal{C}(F_2|H_2) + 1$. Let $F_1^- = (F_1|H_1) \cup (F_2|H_2) \cup \{\{u\}\}$ and $F_2^- = (F_2|H_1) \cup (F_1|H_2) \cup \{\{v\}\}$, see Figure 11. Then $\mathcal{C}(F_1^-) = \mathcal{C}(F_1|H_1) + \mathcal{C}(F_2|H_2) + 1 < \mathcal{C}(F_1|H_1) + \mathcal{C}(F_1|H_2) + 1 = \mathcal{C}(F_1) = k+1$ and $\mathcal{C}(F_2^-) = \mathcal{C}(F_2|H_1) + \mathcal{C}(F_1|H_2) + 1 < \mathcal{C}(F_1|H_1) + \mathcal{C}(F_1|H_2) + 1 = \mathcal{C}(F_1) = k+1$.
If $e \in E(F_2)$, we have $\mathcal{C}(F_1|H_1) \geq \mathcal{C}(F_2|H_1) + 2$ and $\mathcal{C}(F_1|H_2) \geq \mathcal{C}(F_2|H_2) + 2$. Let $F_1^- = (F_1|H_1) \cup (F_2|H_2) \cup \{\{v\}\}$ and $F_2^- = (F_2|H_1) \cup (F_1|H_2) \cup \{\{u\}\}$. Then $\mathcal{C}(F_1^-) = \mathcal{C}(F_1|H_1) + \mathcal{C}(F_2|H_2) + 1 < \mathcal{C}(F_1|H_1) + \mathcal{C}(F_1|H_2) = \mathcal{C}(F_1) = k+1$ and $\mathcal{C}(F_2^-) = \mathcal{C}(F_2|H_1) + \mathcal{C}(F_1|H_2) + 1 < \mathcal{C}(F_1|H_1) + \mathcal{C}(F_1|H_2) = \mathcal{C}(F_1) = k+1$.
In both cases, $\max\{\mathcal{C}(F_1^-), \mathcal{C}(F_2^-)\} \leq k$, thus $H - e \subseteq_I K_k \square K_k$.

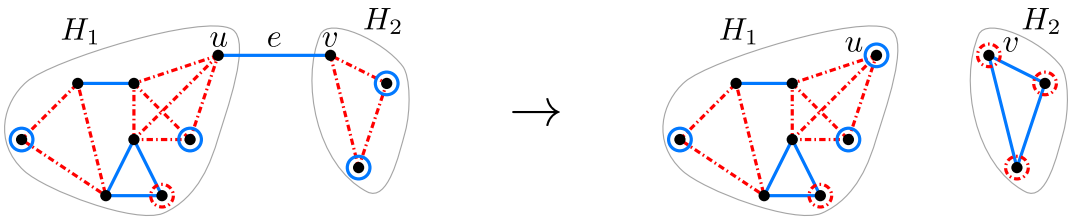

Figure 11: Illustration for the proof of Theorem 12 (Claim 2, case of $H^-$, b1)).

Let b2) hold for $H$ and let $w, T, H_1, H_2$ be as in the definition of class $\mathcal{B}$.
If $T$ is a component of $F_1$, then let $F_1^- = (F_1|H_1) \cup (F_2|H_2) \cup \{\{u,w\}\}$ and $F_2^- = (F_2|H_1) \cup (F_1|H_2) \cup \{\{v,w\}\}$, see Figure 12.
If $T$ is a component of $F_2$, then let $F_1^- = (F_1|H_1) \cup (F_2|H_2) \cup \{\{v,w\}\}$ and $F_2^- = (F_2|H_1) \cup (F_1|H_2) \cup \{\{u,w\}\}$.

In both cases, $\mathcal{C}(F_1^-) = \mathcal{C}(F_1|H_1) + \mathcal{C}(F_2|H_2) + 1 < \mathcal{C}(F_1|H_1) + \mathcal{C}(F_1|H_2) + 1 = \mathcal{C}(F_1) = k+1$ and $\mathcal{C}(F_2^-) = \mathcal{C}(F_2|H_1) + \mathcal{C}(F_1|H_2) + 1 < \mathcal{C}(F_1|H_1) + \mathcal{C}(F_1|H_2) + 1 = \mathcal{C}(F_1) = k+1$. Thus $\max\{\mathcal{C}(F_1^-), \mathcal{C}(F_2^-)\} \leq k$, implying that $H - e \subseteq_I K_k \square K_k$.

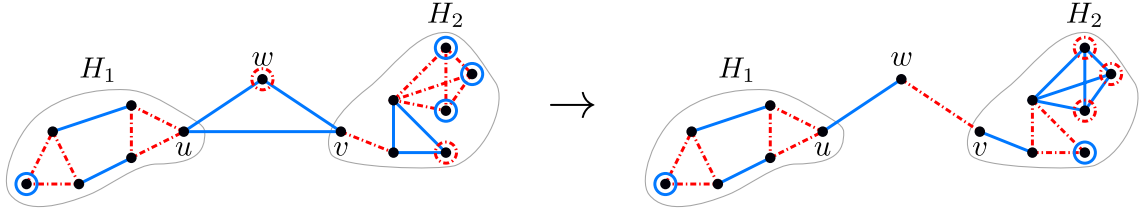

Figure 12: Illustration for the proof of Theorem 12 (Claim 2, case of $H^-$, b2)).

To summarize, we have that $H \not\subseteq_I K_k \square K_k$, $H - e \subseteq_I K_k \square K_k$, and $H + f \subseteq_I K_k \square K_k$. Thus by Lemma 14, $H$ is good. This completes the proof of Claim 2.

**Claim 3.** If $H$ is good then $H$ is either chipped or $H \in \mathcal{A} \cap \mathcal{B}$.

Let $H$ be a good graph. By Lemma 14, there exists $k \in \mathbb{N}$ so that $H \not\subseteq_I K_k \square K_k$ and there are edges $e \in E(H)$ and $f \in E(\overline{H})$ such that $H - e \subseteq_I K_k \square K_k$ and $H + f \subseteq_I K_k \square K_k$. In particular, $H - e$ and $H + f$ are 2-Hamming. Since $H \not\subseteq_I K_k \square K_k$, we have that either $H$ is not 2-Hamming, i.e., $\dim(H) > 2$, or in any 2-Hamming decomposition of $H$, one of the families has at least $k+1$ components, see Lemma 17.

Let $\{F_1^+, F_1^+\}$ be an optimal 2-Hamming decomposition of $H + f$ and $\{F_1^-, F_1^-\}$ be an optimal 2-Hamming decomposition of $H - e$. By Lemma 17, $\max\{\mathcal{C}(F_1^+), \mathcal{C}(F_2^+)\} \leq k$ and $\max\{\mathcal{C}(F_1^-), \mathcal{C}(F_2^-)\} \leq k$.

Let $K(e)$ be the maximal clique of $H$ that contains $e$ and $K(f)$ be the maximal clique of $H + f$ that contains $f$. Recall that by Lemma 16, $H - e$ cannot contain an induced copy of $K_4$ minus an edge. It implies that $|K(e)| \leq 3$, and that $|K(f)| \leq 4$ since otherwise $H$ contains a copy of $K_5$ minus an edge and thus $H - e'$ contains an induced $K_4$ minus an edge for any $e' \in E(H)$.

**Claim 3.1.** If $H$ is good and $\dim(H) > 2$ then $H$ is chipped.

Assume that $H$ is good and $\dim(H) > 2$. We shall show that $H$ is chipped.

Without loss of generality, $f \in E(F_1^+)$ and thus $K = K(f)$ is a component of $F_1^+$.
If $|K| = 2$ then $H$ is 2-Hamming, a contradiction.
If $|K| = 3$, let $K = \{\{u, v, w\}\}$ and $f = uv$. If $\deg(w) \geq 3$, then $H$ is chipped of type ch1). If $\deg(w) = 2$, then $\{\{w\}\}$ is a component of $H - \{uw, wv\}$ and $u, v$ are in the same connected component of $H - \{uw, wv\}$, otherwise $H$ would be 2-Hamming. Thus $H$ is chipped of type ch2).

Recall that by Lemma 16, $H - e$ cannot contain an induced copy $K_{1,3}$ or $K_4$ minus an edge.

Finally let $|K| = 4$. Then $e$ has to be induced by $V(K)$, otherwise $H - e$ contains an induced copy of $K_4$ minus an edge, a contradiction. Hence $V(K)$ induces exactly 4 edges in $H - e$, i.e., $V(K)$ induces either a $C_4$ or a triangle with a pendant edge in $H - e$.

Let $V(K)$ induce $C_4$ in $H - e$. We claim that $H - e$ is isomorphic to $C_4$ then. Otherwise, since $H$ is connected, there is a vertex not in $V(K)$ adjacent to a vertex in $V(K)$. By Lemma 15, it is adjacent to exactly one vertex of $K$ and thus $H - e$ contains an induced copy of $K_{1,3}$, a contradiction. Hence $H - e$ is isomorphic to $C_4$ and thus $H$ is isomorphic to $K_4$ minus an edge, so $H$ is chipped of type ch3).

Let $V(K)$ induce a triangle with a pendant edge in $H - e$. Let $K = \{\{u, v, w, x\}\}$ such that $w$ has degree 3 and $u$ has degree 1 in the graph induced by $V(K)$ in $H - e$. If there is a vertex not in $V(K)$ adjacent to $w$, then by Lemma 15 (applied to $H + f$), it is not adjacent to any vertex in $V(K) \setminus \{w\}$ and thus there is an induced $K_{1,3}$ in $H - e$ centered in $w$, a contradiction. Thus $\deg_H(w) = 3$. Let $\{F_1^-, F_2^-\}$ be a 2-Hamming decomposition of $H - e$. By Lemma 15, $\{\{u, w\}\}$ and $\{\{v, w, x\}\}$ are components of $F_1^-$ or $F_2^-$. Then without loss of generality, $\{\{v, w, x\}\}$ is a component of $F_1^-$ and $\{\{u, w\}\}$ is a component of $F_2^-$. Let $H''$ be the graph obtained from $H$ by deleting the edges induced by $V(K)$ and let $K(u), K(v)$ and $K(x)$ denote the maximal cliques in $H''$ that contain the vertices $u, v$ and $x$ respectively (they might be singletons). Then in $H + f$, $K(u), K(v)$, and $K(x)$ must all belong to $F_2^+$ but in $H - e$, $K(u)$ belongs to $F_1^-$, and both $K(v), K(x)$ belong to $F_2^-$. By Lemma 18, $K(u)$ is not in the same connected component of $H''$ as $K(v)$ or $K(x)$. Therefore the connected component $\widetilde{H}$ of $H''$ that contains $u$ satisfies $V(\widetilde{H}) \cap K = \{u\}$ and thus $H$ is chipped of type ch3).

This proves that $H$ is chipped.

**Claim 3.2. If $\dim(H) = 2$ and $H$ is good then $H \in \mathcal{A} \cap \mathcal{B}$.**

If $\dim(H) = 2$ and $H$ is good, there exists $k \in \mathbb{N}$ so that $H \not\subseteq_I K_k \square K_k$, $H \subseteq_I K_{k+1} \square K_{k+1}$ and there are edges $e \in E(H)$ and $f \in E(\overline{H})$ such that $H - e \subseteq_I K_k \square K_k$ and $H + f \subseteq_I K_k \square K_k$.

By Lemma 17, it means that in any 2-Hamming decomposition of $H$, one of the families has at least $k + 1$ components, but after deleting or adding a suitable edge, we can decompose $H$ with at most $k$ components in each family. Let $\{F_1, F_2\}$ be a 2-Hamming decomposition of $H$.

Recall that $K(e)$ denotes the maximal clique of $H$ that contains $e$ and $K(f)$ denotes the maximal clique of $H + f$ that contains $f$. To avoid redundant arguments for $e$ and $f$, we introduce a variable $g \in \{e, f\}$. Denote $K(g)$ by $K$ and the endpoints of $g$ by $u, v$. Let $H'$ be the graph obtained from $H$ by deleting the edges of $K$. Since $H$ and $H - e$ cannot contain an induced copy of $K_4$ minus an edge, $2 \leq |K| \leq 3$. If $|K| = 3$, then let $K = \{u, v, w\}$ and observe that $\deg_{H'}(w) = 0$, otherwise it follows from Lemma 15, that $H$ or $H - e$ would contain an induced $K_{1,3}$ centered at $w$. Let $K(u)$ and $K(v)$ denote the maximal cliques of $H'$ containing $u$ and $v$ respectively. Observe that $N[u] = V(K \cup K(u))$ and $N[v] = V(K \cup K(v))$, otherwise $H$ would contain an induced $K_{1,3}$. For a 2-Hamming decomposition $\{Q_1, Q_2\}$ of $H$, $H - e$, or $H + f$, let $\{Q'_1, Q'_2\}$ denote the 2-Hamming decomposition of $H' = H - E(K)$ such that $Q'_1$ contains $Q_1|H'$ and $Q'_2$ contains $Q_2|H'$ as subgraphs. Note that $Q'_i$ is obtained from $Q_i|H'$ by possibly adding the singleton vertices $\{\{u\}\}$ or $\{\{v\}\}$.

**Verification of $H \in \mathcal{A}$.**

Let $g = f$. We shall show that $H$ is unbalanced and has property a1) or a2).

If $|K(f)| = 2$, then since the 2-Hamming decomposition of $H$ is unique, the endpoints of $f$ have to be singletons in $F_1$ or $F_2$, say in $F_1$. Since all the components of $\{F_1, F_2\}$, except for the ones containing the endpoints of $f$, are the same as in $\{F_1^+, F_2^+\}$, we have $\mathcal{C}(F_1^+) = \mathcal{C}(F_1) - 1$ and $\mathcal{C}(F_2^+) = \mathcal{C}(F_2)$. We also know that $\max\{\mathcal{C}(F_1), \mathcal{C}(F_2)\} = k + 1$ and $\max\{\mathcal{C}(F_1^+), \mathcal{C}(F_2^+)\} = k$, therefore it follows that $\mathcal{C}(F_2) \leq k$ and $\mathcal{C}(F_1) = k + 1$. Thus $H$ is unbalanced, has property a1) and so $H \in \mathcal{A}$.

If $|K(f)| = 3$, then without loss of generality, $K(f)$ is a component of $F_1^+$. Recall that in this case $\deg_{H'}(w) = 0$, i.e., $\deg_H(w) = 2$. If follows that $K(u)$ and $K(v)$ are components in $(F_2^+)'$, but on the other hand, they belong to different families in the 2-Hamming decomposition $\{F'_1, F'_2\}$ of $H'$. By Lemma 18, $K(u)$ and $K(v)$ are in different components of $H'$. Let $H_1, H_2, \{\{w\}\}$ be the components of $H'$ with $u \in V(H_1)$ and $v \in V(H_2)$. By the uniqueness of the 2-Hamming decomposition of connected graphs, the only 2-Hamming decomposition of $H$ is $\{\widetilde{F}_1, \widetilde{F}_2\}$, where $\widetilde{F}_1 = (F_1^+|H_1) \cup (F_2^+|H_2) \cup \{\{u, w\}\}$ and $\widetilde{F}_2 = (F_2^+|H_1) \cup (F_1^+|H_2) \cup \{\{v, w\}\}$.

Let $n_i = \mathcal{C}(F_i^+|H_1)$ and $l_i = \mathcal{C}(F_i^+|H_2)$. Then

$$\mathcal{C}(\widetilde{F}_i|H_1) = \mathcal{C}(F_i^+|H_1) = n_i \quad \text{and} \quad \mathcal{C}(\widetilde{F}_i|H_2) = \mathcal{C}(F_{3-i}^+|H_2) = l_{3-i}. \tag{1}$$

We see that $\mathcal{C}(F_i^+) = n_i + l_i + 1$ and $\mathcal{C}(\widetilde{F}_i) = n_i + l_{3-i} + 1$.
By Lemma 17, $\max\{\mathcal{C}(F_1^+), \mathcal{C}(F_2^+)\} \leq k$ and $\max\{\mathcal{C}(\widetilde{F}_1), \mathcal{C}(\widetilde{F}_2)\} = k + 1$, which gives

$$\max\{n_1 + l_1 + 1, n_2 + l_2 + 1\} < \max\{n_1 + l_2 + 1, n_2 + l_1 + 1\}. \tag{2}$$

Let $i \in \{1, 2\}$ such that $\mathcal{C}(\widetilde{F}_i) \geq \mathcal{C}(\widetilde{F}_{3-i})$, i.e., $n_i + l_{3-i} + 1 \geq n_{3-i} + l_i + 1$. Then (2) imples that $\max\{n_1 + l_1 + 1, n_2 + l_2 + 1\} < n_i + l_{3-i} + 1$ and thus $n_i > n_{3-i}$ and $l_{3-i} > l_i$. Let $F_1 = \widetilde{F}_i$ and $F_2 = \widetilde{F}_{3-i}$. Then by (1), $n_i > n_{3-i}$ means that $\mathcal{C}(F_1|H_1) > \mathcal{C}(F_2|H_1)$ and $l_{3-i} > l_i$ means that $\mathcal{C}(F_1|H_2) > \mathcal{C}(F_2|H_2)$. Finally, we have that $\mathcal{C}(F_1) = \mathcal{C}(F_1|H_1) + \mathcal{C}(F_1|H_2) + 1 > \mathcal{C}(F_2|H_1) + \mathcal{C}(F_2|H_1) + 1 = \mathcal{C}(F_2)$. Thus $H$ is unbalanced, has property a2) and so $H \in \mathcal{A}$.

**Verification of $H \in \mathcal{B}$.**

Let $g = e$. We shall show that $H$ is unbalanced and has property b1) or b2).

Without loss of generality, $K(u)$ is a component of $F_1^-$. It follows that $K(v)$ has to be a component of $F_2^-$, indeed, if $|K(e)| = 3$, it is clear and if $|K(e)| = 2$ and $K(v)$ is a component of $F_1^-$, then $F_1^-$ and $(F_2^- - \{\{u\}\} - \{\{v\}\}) \cup \{\{u, v\}\}$ would form a 2-Hamming decomposition of $H$ with at most $k$ components in both families, a contradiction. On the other hand, Lemma 15 implies that $K(e)$ is a component in one of $\{F_1, F_2\}$, therefore $K(u)$ and $K(v)$ belong the same family in $\{F_1', F_2'\}$. By Lemma 18, $K(u)$ and $K(v)$ are in different connected components of $H'$. Let $H_1, H_2$ denote the components of $H'$ with $u \in V(H_1)$ and $v \in V(H_2)$.

As observed before, $N[u] = V(K(e) \cup K(u))$ and $N[v] = V(K(e) \cup K(v))$. Thus if $|K(e)| = 2$, then $\{\{u\}\}$ and $\{\{v\}\}$ are singleton components of $F_2^-$ and $F_1^-$ respectively.

By the uniqueness of the 2-Hamming decomposition of connected graphs, the unique 2-Hamming decomposition $\{\widetilde{F}_1, \widetilde{F}_2\}$ of $H$ is

$$\widetilde{F}_1 = \begin{cases} (F_1^-|H_1) \cup (F_2^-|H_2) & \text{if } |K(e)| = 2, \\ (F_1^-|H_1) \cup (F_2^-|H_2) \cup \{\{w\}\} & \text{if } |K(e)| = 3, \end{cases}$$

$$\widetilde{F}_2 = \begin{cases} ((F_2^-|H_1) - \{\{u\}\}) \cup ((F_1^-|H_2) - \{\{v\}\}) \cup \{\{u,v\}\} & \text{if } |K(e)| = 2, \\ (F_2^-|H_1) \cup (F_1^-|H_2) \cup \{\{u,v,w\}\} & \text{if } |K(e)| = 3. \end{cases}$$

Let $n_i = \mathcal{C}(F_i^-|H_1)$ and $l_i = \mathcal{C}(F_i^-|H_2)$, then

$$\mathcal{C}(\widetilde{F}_1|H_1) = n_1, \qquad \mathcal{C}(\widetilde{F}_1|H_2) = l_2, \tag{3}$$
$$\mathcal{C}(\widetilde{F}_2|H_1) = n_2 - 1, \qquad \mathcal{C}(\widetilde{F}_2|H_2) = l_1 - 1 \qquad \text{if } |K(e)| = 2, \tag{4}$$
$$\mathcal{C}(\widetilde{F}_2|H_1) = n_2, \qquad \mathcal{C}(\widetilde{F}_2|H_2) = l_1 \qquad \text{if } |K(e)| = 3. \tag{5}$$

It follows that for $i \in \{1, 2\}$,

$$\mathcal{C}(F_i^-) = n_i + l_i, \qquad \mathcal{C}(\widetilde{F}_1) = n_1 + l_2, \quad \text{and} \quad \mathcal{C}(\widetilde{F}_2) = n_2 + l_1 - 1 \text{ if } |K(e)| = 2, \tag{6}$$
$$\mathcal{C}(F_i^-) = n_i + l_i + 1, \quad \mathcal{C}(\widetilde{F}_1) = n_1 + l_2 + 1, \quad \text{and} \quad \mathcal{C}(\widetilde{F}_2) = n_2 + l_1 + 1 \text{ if } |K(e)| = 3. \tag{7}$$

By Lemma 17, $\max\{\mathcal{C}(F_1^-), \mathcal{C}(F_2^-)\} \leq k$ and $\max\{\mathcal{C}(\widetilde{F}_1), \mathcal{C}(\widetilde{F}_2)\} = k + 1$, which gives

$$\max\{n_1 + l_1, n_2 + l_2\} < \max\{n_1 + l_2, n_2 + l_1 - 1\} \qquad \text{if } |K(e)| = 2, \tag{8}$$
$$\max\{n_1 + l_1 + 1, n_2 + l_2 + 1\} < \max\{n_1 + l_2 + 1, n_2 + l_1 + 1\} \qquad \text{if } |K(e)| = 3. \tag{9}$$

$$\text{Let } i \in \{1, 2\} \text{ such that } \mathcal{C}(\widetilde{F}_{3-i}) \leq \mathcal{C}(\widetilde{F}_i) \text{ and let } F_1 := \widetilde{F}_i, \ F_2 := \widetilde{F}_{3-i}. \tag{10}$$

Assume that $|K(e)| = 2$. By the definition of $\{\widetilde{F}_1, \widetilde{F}_2\}$, $e \in E(\widetilde{F}_2)$. We shall show that $H$ is unbalanced and that b1) holds, i.e., that if $e \in E(F_j)$, then $\mathcal{C}(F_2|H_1) + j \leq \mathcal{C}(F_1|H_1)$ and $\mathcal{C}(F_2|H_2) + j \leq \mathcal{C}(F_1|H_2)$. It follows from (10) and the fact that $e \in E(\widetilde{F}_2)$ that $e \in E(F_j)$ if and only if $F_j = \widetilde{F}_2$ and $F_{3-j} = \widetilde{F}_1$.

If $j = 1$, then

$$F_1 = \widetilde{F}_2, \ F_2 = \widetilde{F}_1, \text{ and} \tag{11}$$

$$n_1 + l_2 \stackrel{(6)}{=} \mathcal{C}(\widetilde{F}_1) \stackrel{(11)}{=} \mathcal{C}(F_2) \stackrel{(10)}{\leq} \mathcal{C}(F_1) \stackrel{(11)}{=} \mathcal{C}(\widetilde{F}_2) \stackrel{(6)}{=} n_2 + l_1 - 1. \tag{12}$$

From (8) and (12), we obtain $\max\{n_1 + l_1, n_2 + l_2\} < n_2 + l_1 - 1$, which implies $n_1 < n_2 - 1$ and $l_2 < l_1 - 1$. Then $H$ is unbalanced since

$$\mathcal{C}(F_2) \stackrel{(11)}{=} \mathcal{C}(\widetilde{F}_1) \stackrel{(6)}{=} n_1 + l_2 < n_2 + l_2 - 1 < n_2 + l_1 - 1 \stackrel{(6)}{=} \mathcal{C}(\widetilde{F}_2) \stackrel{(11)}{=} \mathcal{C}(F_1).$$

Moreover,

$$\mathcal{C}(F_2|H_1) + 1 \stackrel{(11)}{=} \mathcal{C}(\widetilde{F}_1|H_1) + 1 \stackrel{(3)}{=} n_1 + 1 \leq n_2 - 1 \stackrel{(4)}{=} \mathcal{C}(\widetilde{F}_2|H_1) \stackrel{(11)}{=} \mathcal{C}(F_1|H_1),$$

$$\mathcal{C}(F_2|H_2) + 1 \stackrel{(11)}{=} \mathcal{C}(\widetilde{F}_1|H_2) + 1 \stackrel{(3)}{=} l_2 + 1 \leq l_1 - 1 \stackrel{(4)}{=} \mathcal{C}(\widetilde{F}_2|H_2) \stackrel{(11)}{=} \mathcal{C}(F_1|H_2).$$

Thus $H$ has property b1).

If $j = 2$, then

$$F_1 = \widetilde{F}_1,\ F_2 = \widetilde{F}_2,\ \text{and} \tag{13}$$

$$n_2 + l_1 - 1 \stackrel{(6)}{=} \mathcal{C}(\widetilde{F}_2) \stackrel{(13)}{=} \mathcal{C}(F_2) \stackrel{(10)}{\leq} \mathcal{C}(F_1) \stackrel{(13)}{=} \mathcal{C}(\widetilde{F}_1) \stackrel{(6)}{=} n_1 + l_2. \tag{14}$$

From (8) and (14), we obtain $\max\{n_1 + l_1, n_2 + l_2\} < n_1 + l_2$, which implies $n_2 < n_1$ and $l_1 < l_2$. Then $H$ is unbalanced since

$$\mathcal{C}(F_2) \stackrel{(13)}{=} \mathcal{C}(\widetilde{F}_2) \stackrel{(6)}{=} n_2 + l_1 - 1 < n_2 + l_1 < n_1 + l_2 \stackrel{(6)}{=} \mathcal{C}(\widetilde{F}_1) \stackrel{(13)}{=} \mathcal{C}(F_1).$$

Moreover,

$$\mathcal{C}(F_2|H_1) + 2 \stackrel{(13)}{=} \mathcal{C}(\widetilde{F}_2|H_1) + 2 \stackrel{(4)}{=} n_2 + 1 \leq n_1 \stackrel{(3)}{=} \mathcal{C}(\widetilde{F}_1|H_1) \stackrel{(13)}{=} \mathcal{C}(F_1|H_1),$$

$$\mathcal{C}(F_2|H_2) + 2 \stackrel{(13)}{=} \mathcal{C}(\widetilde{F}_2|H_2) + 2 \stackrel{(4)}{=} l_1 + 1 \leq l_2 \stackrel{(3)}{=} \mathcal{C}(\widetilde{F}_1|H_2) \stackrel{(13)}{=} \mathcal{C}(F_1|H_2).$$

Thus $H$ has property b1).

Assume that $|K(e)| = 3$. We shall show that $H$ is unbalanced and that b2) holds, i.e. that $\mathcal{C}(F_2|H_1) < \mathcal{C}(F_1|H_1)$, $\mathcal{C}(F_2|H_2) < \mathcal{C}(F_1|H_2)$ and $H_3 = \{\{w\}\}$. We showed in the beginning of the proof of Claim 3.2, that $w$ has degree 0 in $H' = H - K(e)$, thus $\{\{w\}\}$ is a component of $H - K(e)$. By the choice of $i$,

$$n_{3-i} + l_i + 1 \stackrel{(7)}{=} \mathcal{C}(\widetilde{F}_{3-i}) \stackrel{(10)}{\leq} \mathcal{C}(\widetilde{F}_i) \stackrel{(7)}{=} n_i + l_{3-i} + 1. \tag{15}$$

From (9) and (15), we get $\max\{n_1 + l_1 + 1, n_2 + l_2 + 1\} < n_i + l_{3-i} + 1$, which implies that $n_i > n_{3-i}$ and $l_{3-i} > l_i$. Then $H$ is unbalanced since

$$\mathcal{C}(F_2) \stackrel{(10)}{=} \mathcal{C}(\widetilde{F}_{3-i}) \stackrel{(7)}{=} n_{3-i} + l_i + 1 < n_i + l_{3-i} + 1 \stackrel{(7)}{=} \mathcal{C}(\widetilde{F}_i) \stackrel{(10)}{=} \mathcal{C}(F_1).$$

Moreover,

$$\mathcal{C}(F_2|H_1) \stackrel{(10)}{=} \mathcal{C}(\widetilde{F}_{3-i}|H_1) \stackrel{(3),(5)}{=} n_{3-i} < n_i \stackrel{(3),(5)}{=} \mathcal{C}(\widetilde{F}_i|H_1) \stackrel{(10)}{=} \mathcal{C}(F_1|H_1),$$

$$\mathcal{C}(F_2|H_2) \stackrel{(10)}{=} \mathcal{C}(\widetilde{F}_{3-i}|H_2) \stackrel{(3),(5)}{=} l_i < l_{3-i} \stackrel{(3),(5)}{=} \mathcal{C}(\widetilde{F}_i|H_2) \stackrel{(10)}{=} \mathcal{C}(F_1|H_2).$$

Thus $H$ has property b2).

Therefore, $H$ is unbalanced, b1) holds if $|K(e)| = 2$, b2) holds if $|K(e)| = 3$ and thus $H \in \mathcal{B}$. This concludes the proof of Claim 3.2.

Now, to prove the main statement of the theorem, assume first that $H$ is good. If $\dim(H) > 2$, then by Claim 3.1 $H$ is chipped. If $\dim(H) \leq 2$, then by Claim 3.2 $H \in \mathcal{A} \cap \mathcal{B}$. The other way around, assume that $H$ is chipped or $H \in \mathcal{A} \cap \mathcal{B}$ then by Claims 1 and 2 $H$ is good. □

**Proof of Proposition 4.** We shall show that $H$ is either in the class $\mathcal{T}(n)$ or chipped and so by Corollary 2 and Theorem 12, there is a graph $G$ so that $G \hookrightarrow\!\!\!\!\!\rightarrow H$. If $H$ is a spider with at least 4 legs, then $H \in \mathcal{T}(n)$. If $H$ is a spider with 3 legs, then $H$ is chipped of type ch1). If $H \cong \widetilde{C}_{2k+1}$, then $H$ is chipped of type ch2). Finally, if $H$ is a union of a cycle $C$ and a path $P$ such that $P$ has length at least 1 and shares exactly one vertex with $C$ which is an endpoint of it, then $H$ is chipped of type ch1). □

# 6  Further results

For a family $\mathcal{F}$ of graphs, we use $\mathrm{Forb}(\mathcal{F})$ to denote a set of all graphs that do not contain an induced copy of any graph from $\mathcal{F}$. When $\mathcal{F} = \{H\}$, we write $\mathrm{Forb}(\{H\}) = \mathrm{Forb}(H)$. We say that $G$ is $\mathcal{F}$-*induced-saturated* and write $G \hookrightarrow\!\!\!\!\!\rightarrow \mathcal{F}$ if $G \in \mathrm{Forb}(\mathcal{F})$, but deleting any edge from $G$ or adding any edge from $\overline{G}$ to $G$ results in an induced copy of a member of $\mathcal{F}$.

Using the structural descriptions of the classes $\mathrm{Forb}(C_4, \overline{C}_4, P_4)$ (due to Chvátal and Hammer [2]) and $\mathrm{Forb}(C_4, \overline{C}_4, C_5)$ (due to Földes and Hammer [6]), Behrens et al.[1] showed that there is no graph $G$ so that $G \hookrightarrow\!\!\!\!\!\rightarrow \{\overline{C}_4, C_4, P_4\}$ or $G \hookrightarrow\!\!\!\!\!\rightarrow \{\overline{C}_4, C_4, C_5\}$. The later example shows that even though there are graphs $G_1, G_2, G_3$ so that $G_1 \hookrightarrow\!\!\!\!\!\rightarrow \overline{C}_4$, $G_2 \hookrightarrow\!\!\!\!\!\rightarrow C_4$ and $G_3 \hookrightarrow\!\!\!\!\!\rightarrow C_5$, there is no induced-saturated graph for the family $\{\overline{C}_4, C_4, C_5\}$.

In this section, we give a short proof of the fact, first proved by Martin and Smith [9], that there is no $P_4$-induced-saturated graph and we further extend the list of families for which there is no induced-saturated graph, see Theorem 19. We generalize the ideas used in these proofs and provide some saturation properties of prime graphs, see Theorem 20.

First we need a few definitions. We say that a vertex of a graph is *full* if it is adjacent to all other vertices of the graph. A *complement reducible graph*, shortly *cograph*, is defined recursively from a single vertex graph by taking disjoint unions and complements, i.e., $K_1$ is a cograph and for cographs $G_1, G_2$, the graphs $G_1 \dot\cup G_2$, $\overline{G}_1$, $\overline{G}_2$ are also cographs. Corneil et al.[3] showed that a graph $G$ is a cograph if and only if $G \in \mathrm{Forb}(P_4)$. A *homogeneous set* of an $n$-vertex graph $G$ is a set $X \subsetneq V(G)$ with $2 \leq |X| \leq n-1$, such that every vertex $v \in V(G) \setminus X$ is either complete to $X$, i.e., $X \subseteq N(v)$ or anticomplete to $X$, i.e., $X \cap N(v) = \emptyset$. We say that a graph $H$ is *prime* if it does not have a homogeneous set, that is, for any proper subset $X$ of vertices with $|X| \geq 2$, there is a vertex $v \notin X$ such that $v$ sends both an edge and a non-edge to $X$. Note that if $H$ is a prime graph on at least 2 vertices, then $H$ and $\overline{H}$ are connected, in particular $H$ does not have an isolated or a full vertex.

**Theorem 19.** *There is no graph $G$ such that $G \hookrightarrow\!\!\!\!\!\rightarrow P_4$ and there is no graph $G$ such that $G \hookrightarrow\!\!\!\!\!\rightarrow \{P_4, C_4\}$.*

*Proof.* Assume that $G$ is a graph on the smallest number of vertices such that $G \hookrightarrow\!\!\!\!\!\rightarrow P_4$. Then $|G| \geq 4$ and since $P_4 \not\subseteq_I G$, $G$ is a cograph. In particular, $G$ or its complement is disconnected. Assume without loss of generality that $G$ is disconnected and let $G'$ be a union of some connected components of $G$ such that $G' \neq G$ and $|G'| \geq 2$. Deleting an edge of $G'$ or adding an edge from the complement of $G'$ creates a copy of $P_4$ in $G$. Since $P_4$ is connected, this copy cannot contain vertices from both $V(G')$ and $V(G) \setminus V(G')$, i.e., it lies completely in $G'$. Thus $G' \hookrightarrow\!\!\!\!\!\rightarrow P_4$, which contradicts the minimality of $G$.

Assume that $G = (V, E)$ is a graph on the smallest number of vertices so that $G \hookrightarrow\!\!\!\!\!\rightarrow \{P_4, C_4\}$. Then $G$ is $\{P_4, C_4\}$-free and by the minimality of $G$, $G$ is connected. Wolk [12] proved that every finite, connected $\{P_4, C_4\}$-free graph has a full vertex. Thus $G$ has a full vertex, denote it by $v$ and let $G' = G - v$. Since neither $P_4$ nor $C_4$ contains a full vertex, deleting an edge of $G'$ or adding an edge to $G'$ from its complement results in an induced copy of $P_4$ or $C_4$ that does not use $v$, i.e. it lies in $G'$. Therefore $G' \hookrightarrow\!\!\!\!\!\rightarrow \{P_4, C_4\}$ which contradicts the minimality of $G$. □

For two graphs $G_1 = (V_1, E_1)$ and $G_2 = (V_2, E_2)$, the *blowup product* of $G_1$ and $G_2$, denoted by $G_1 \star G_2$ is a graph on vertex set $V = V_1 \times V_2$ where two vertices $(u_1, v_1), (u_2, v_2) \in V$ are adjacent if either $u_1 u_2 \in E_1$ or $u_1 = u_2$ and $v_1 v_2 \in E_2$. Informally, $G_1 \star G_2$ is obtained by blowing up $G_1$ and inserting a copy of $G_2$ in each blob.

**Theorem 20.** *Let $H$ be a prime graph for which there exists an $H$-induced-saturated graph $G$. Then $G \star G \hookrightarrow\!\!\!\!\!\rightarrow H$ and there is a prime graph $G'$ so that $G' \hookrightarrow\!\!\!\!\!\rightarrow H$.*

*Proof.* Let $G \hookrightarrow\!\!\!\!\!\rightarrow H$.

We shall show first that $G \star G \hookrightarrow\!\!\!\!\!\rightarrow H$. We need to verify that $H \not\subseteq_I G \star G$ and that adding or deleting any edge results in an induced copy of $H$. Let $X \subseteq V(G \star G)$ of size $|H|$. If $X$ intersects each blob of $G \star G$ in at most one vertex or $X$ is contained in one blob of $G \star G$, then $X$ induces a subgraph of $G$ and thus $(G \star G)[X]$ is not isomorphic to $H$. If there exists a blob $B$ so that $X \not\subseteq V(B)$ and $X \cap B = Y$ with $|Y| \geq 2$, then $X$ cannot induce a prime graph since a vertex of $X$ that is not in $B$ is either complete or anticomplete to $Y$. If we add or delete an edge between vertices that are contained in the same blob of $G \star G$, then since the blob induces $G$ and $G \hookrightarrow\!\!\!\!\!\rightarrow H$, this change results in an induced copy of $H$ within the blob. If we add or delete an edge between vertices that are contained in different blobs of $G \star G$, then this pair belongs to a copy $\widetilde{G}$ of $G$ obtained by taking a single vertex from each blob. Thus this change corresponds to an edge addition or deletion in $\widetilde{G}$ that results in an induced copy of $H$.

Next we shall show that there is a prime graph $G'$ so that $G' \hookrightarrow\!\!\!\!\!\rightarrow H$. Suppose that $G$ is not a prime graph and consider a minimal homogeneous set $S$ in $G$. Let $G' = G[S]$.

Observe that if $X \subsetneq V(G')$ is a homogeneous set in $G'$ then it is also a homogeneous set in $G$. By the minimality of $S$, we conclude that $G'$ does not have a homogeneous set, i.e., $G'$ is a prime graph. Since $G \hookrightarrow\!\!\!\!\!\!\twoheadrightarrow H$, deleting an edge in $G'$ or adding an edge from the complement of $G'$ creates an induced copy $H'$ of $H$. Since $H$ is a prime graph, $H'$ cannot contain vertices from both $S$ and $V(G) \setminus S$ therefore $H' \subseteq G'$ and hence $G' \hookrightarrow\!\!\!\!\!\!\twoheadrightarrow H$. $\square$

**Corollary 21.** *If $H$ is a prime graph such that each $G \in \mathrm{Forb}(H)$ is nonprime, then there is no $H$-induced-saturated graph. In particular there is no $P_4$-induced-saturated graph.*

**Corollary 22.** *If $H$ is a prime graph and there is a graph $G$ such that $G \hookrightarrow\!\!\!\!\!\!\twoheadrightarrow H$, then there is an arbitrarily large graph $G'$ so that $G' \hookrightarrow\!\!\!\!\!\!\twoheadrightarrow H$.*

# 7 Concluding remarks

We have shown that graphs from certain families of graphs have Hamming graphs as induced saturating graphs. Our families do not include paths, hence it remains open whether there is an induced saturating graph for any path of length at least 4.

Note that the blowup product results in a nonprime graph. It is interesting to ask whether for all prime graphs $H$ for which there is an $H$-induced-saturated graph, there exists an arbitrarily large prime graph $G$ so that $G \hookrightarrow\!\!\!\!\!\!\twoheadrightarrow H$.